\vsize=9.0in\voffset=1cm
\looseness=2


\message{fonts,}

\font\tenrm=cmr10
\font\ninerm=cmr9
\font\eightrm=cmr8
\font\teni=cmmi10
\font\ninei=cmmi9
\font\eighti=cmmi8
\font\ninesy=cmsy9
\font\tensy=cmsy10
\font\eightsy=cmsy8
\font\tenbf=cmbx10
\font\ninebf=cmbx9
\font\tentt=cmtt10
\font\ninett=cmtt9

\font\ninesl=cmsl9
\font\eightsl=cmsl8

\font\nineit=cmti9
\font\eightit=cmti8

\skewchar\ninei='177 \skewchar\eighti='177
\skewchar\ninesy='60 \skewchar\eightsy='60

\def\eightpoint{\def\rm{\fam0\eightrm} 
\normalbaselineskip=9pt
\normallineskiplimit=-1pt
\normallineskip=0pt

\textfont0=\eightrm \scriptfont0=\sevenrm \scriptscriptfont0=\fiverm
\textfont1=\ninei \scriptfont1=\seveni \scriptscriptfont1=\fivei
\textfont2=\ninesy \scriptfont2=\sevensy \scriptscriptfont2=\fivesy
\textfont3=\tenex \scriptfont3=\tenex \scriptscriptfont3=\tenex
\textfont\itfam=\eightit  \def\it{\fam\itfam\eightit} 
\textfont\slfam=\eightsl \def\sl{\fam\slfam\eightsl} 

\setbox\strutbox=\hbox{\vrule height6pt depth2pt width0pt}%
\normalbaselines \rm}

\def\ninepoint{\def\rm{\fam0\ninerm} 
\textfont0=\ninerm \scriptfont0=\sevenrm \scriptscriptfont0=\fiverm
\textfont1=\ninei \scriptfont1=\seveni \scriptscriptfont1=\fivei
\textfont2=\ninesy \scriptfont2=\sevensy \scriptscriptfont2=\fivesy
\textfont3=\tenex \scriptfont3=\tenex \scriptscriptfont3=\tenex
\textfont\itfam=\nineit  \def\it{\fam\itfam\nineit} 
\textfont\slfam=\ninesl \def\sl{\fam\slfam\ninesl} 
\textfont\bffam=\ninebf \scriptfont\bffam=\sevenbf
\scriptscriptfont\bffam=\fivebf \def\bf{\fam\bffam\ninebf} 
\textfont\ttfam=\ninett \def\tt{\fam\ttfam\ninett} 

\normalbaselineskip=11pt
\setbox\strutbox=\hbox{\vrule height8pt depth3pt width0pt}%
\let \smc=\sevenrm \let\big=\ninebig \normalbaselines
\parindent=1em
\rm}

\def\tenpoint{\def\rm{\fam0\tenrm} 
\textfont0=\tenrm \scriptfont0=\ninerm \scriptscriptfont0=\fiverm
\textfont1=\teni \scriptfont1=\seveni \scriptscriptfont1=\fivei
\textfont2=\tensy \scriptfont2=\sevensy \scriptscriptfont2=\fivesy
\textfont3=\tenex \scriptfont3=\tenex \scriptscriptfont3=\tenex
\textfont\itfam=\nineit  \def\it{\fam\itfam\nineit} 
\textfont\slfam=\ninesl \def\sl{\fam\slfam\ninesl} 
\textfont\bffam=\ninebf \scriptfont\bffam=\sevenbf
\scriptscriptfont\bffam=\fivebf \def\bf{\fam\bffam\tenbf} 
\textfont\ttfam=\tentt \def\tt{\fam\ttfam\tentt} 

\normalbaselineskip=11pt
\setbox\strutbox=\hbox{\vrule height8pt depth3pt width0pt}%
\let \smc=\sevenrm \let\big=\ninebig \normalbaselines
\parindent=1em
\rm}

\message{fin format jgr}
\magnification=1200
\font\Bbb=msbm10
\def\C{\hbox{\Bbb C}}
\def\R{\hbox{\Bbb R}}
\def\pa{\partial}
\def\b{\backslash}
\def\ep{\varepsilon}
\def\v{\varphi}

\vskip 2 mm
\centerline{\bf On the reconstruction of conductivity}
\centerline{\bf of bordered two-dimensional surface in $\R^3$}
\centerline{\bf from electrical currents measurements on its boundary}
\vskip 2 mm
\centerline{\bf G.M. Henkin$^1$ and R.G. Novikov$^2$}
\vskip 4 mm

\noindent
$^1${\ninerm Universit\'e Pierre et Marie Curie, case 247, 4 place Jussieu,
75252, Paris, France}

\noindent
{\ninerm e-mail: henkin@math.jussieu.fr}

\vskip 2 mm

\noindent
$^2${\ninerm CNRS (UMR 7641), Centre de Math\'ematiques Appliqu\'ees,
Ecole Polytechnique,}

\noindent
{\ninerm 91128, Palaiseau, France}

\noindent
{\ninerm e-mail: novikov@cmap.polytechnique.fr}
\vskip 4 mm
{\bf Abstract.}

An electrical potential $U$ on a bordered real surface $X$ in $\R^3$ with
isotropic conductivity function $\sigma>0$ satisfies equation
$d(\sigma d^cU)\big|_X=0$, where $d^c=i(\bar\pa-\pa)$,
$d=\bar\pa+\pa$ are real operators
associated with complex (conforme) structure on $X$  induced by Euclidien
metric of $\R^3$. This paper gives exact reconstruction of conductivity
function $\sigma$ on $X$ from Dirichlet-to-Neumann mapping
$U\big|_{bX}\to\sigma d^c U\big|_{bX}$. This paper extends to the case of
the Riemann surfaces the reconstruction schemes of R.Novikov [N2] and of
A.Bukhgeim [B], given for the case $X\subset\R^2$. The paper extends and
 corrects the statements of [HM], where the inverse
boundary value problem on
the Riemann surfaces was firstly considered.


\vskip 2 mm
Mathematics Subject Classification (2000) 32S65, 32V20, 35R05, 35R30, 58J32,

\noindent
81U40.

\vskip 2 mm
{\bf 0. Introduction}

{\sl
0.1. Reduction of inverse boundary value problem on a surface in $\R^3$ to
the corresponding problem on affine algebraic Riemann surface in $\C^3$. }

Let $X$ be bordered oriented two-dimensional manifold in $\R^3$. Manifold $X$
is equiped by complex (conformal) structure induced by Euclidean metric
of $\R^3$. We say that $X$ possesses an isotropic conductivity function
$\sigma>0$, if any electric potential $u$ on $bX$ generates electrical
potential $U$ on $X$, solving the Dirichlet problem:
$$U\big|_{bX}=u\ \ {\rm and}\ \ d\sigma d^cU\big|_X=0,\eqno(0.1)$$
where $d^c=i(\bar\pa-\pa)$, $d=\bar\pa+\pa$ and the Cauchy-Riemann operator
$\bar\pa$ corresponds to complex (conformal) structure on $X$. Inverse
conductivity problem consists in the  reconstruction of $\sigma\big|_X$
from the mapping potential
$U\big|_{bX}\to\ {\rm current}\ j=\sigma d^cU\big|_{bX}$ for solutions of
(0.1). This mapping is called Dirichlet-to-Neumann mapping.

This problem is the special case of the following more general inverse
boundary value problem, going back to I.M.Gelfand [Ge] and A.Calderon [C]:
to find potential (2-forme) $q$ on $X$ in the equation
$$d d^c\psi=q\psi \eqno(0.2)$$
from knowledge of Dirichlet-to-Neumann mapping
$\psi\big|_{bX}\to  d^c\psi\big|_{bX}$ for solutions of (0.2). Equation (0.2)
is called in some context by stationary  Schr\"odinger equation,
in other context by monochromatic acoustic equation etc.
Equation (0.1) can be reduced to the equation (0.2) with

\noindent
$q={dd^c\sqrt{\sigma}\over \sqrt{\sigma}}$ by the
substitution $\psi=\sqrt{\sigma}U$.

Let restriction of Euclidean metric of $\R^3$ on $X$ have (in local
coordinates) the form
$$ds^2=Edx^2+2Fdx dy+Gdy^2=Adz^2+2Bdzd\bar z +\bar Ad\bar z^2,$$
where $z=x+iy$, $B={{E+G}\over 4}$, $A={{E-G-2iF}\over 4}$. Put
$\mu={\bar A\over {B+\sqrt{B^2-|A|^2}}}$.
By classical results (going back to Gauss and Riemann) one can construct
holomorphic embedding $\v:\ X\to\C^3$, using some solution of Beltrami
equation: $\bar\pa\v=\mu\pa\v$ on $X$. Moreover, embedding $\v$ can be chosen
in such a way that $\v(X)$ belongs to smooth algebraic curve $V$ in $\C^3$.
Using
existence of embedding $\v$ we can identify further $X$ with $\v(X)$.

\vskip 2 mm
{\sl
0.2. Reconstruction schemes for the case $X\subset\R^2\simeq\C$. }

For the case $X=\Omega\subset\R^2$ the exact reconstruction scheme for
formulated inverse problems was given in [N2], [N3] under some restriction
(smallness assumption) for $\sigma$ or $q$ (see Corollary 2 of [N2])
. For the case of inverse
conductivity problem, see (0.1), (0.2), when $q={dd^c\sqrt{\sigma}\over \sqrt{\sigma}}$,
restriction on $\sigma$ in this scheme was eliminated by A.Nachman [Na] by
the reduction to the equivalent question for the first order system studied by
 R.Beals and R.Coifman [BC2].
Recently A.Bukhgeim [B] has found new original reconstruction scheme for
inverse boundary value problem, see  (0.2), without smallness assumption on
$q$.

In a particular case, the scheme of [N2] for the inverse conductivity problem
 consists in the following. Let $\sigma(x)>0$ for $x\in\bar\Omega$ and
$\sigma\in C^{(2)}(\bar\Omega)$. Put $\sigma(x)=1$ for $x\in\R^2\b\bar\Omega$.

Let $q={dd^c\sqrt{\sigma}\over \sqrt{\sigma}}$.

From L.Faddeev [F1] result it follows:\ $\exists$ compact set $E\subset\C$
 such that for each

\noindent
$\lambda\in\C\b E$ there exists a unique solution
$\psi(z,\lambda)$ of the equation
$dd^c\psi=q\psi={dd^c\sqrt{\sigma}\over \sqrt{\sigma}}\psi$, with asymptotics
$$\psi(z,\lambda)e^{-\lambda z}\buildrel \rm def \over =\mu(z,\lambda)=
1+o(1),\ z\to\infty.$$
Such solution can be found from the integral equation
$$\mu(z,\lambda)=1+{i\over 2}\int\limits_{\xi\in\Omega}g(z-\xi,\lambda)
{\mu(\xi,\lambda)dd^c\sqrt{\sigma}\over \sqrt{\sigma}},\eqno(0.3)$$
where the function
$$g(z,\lambda)={i\over (2\pi)^2}\int\limits_{w\in\C}
{e^{\lambda w-\bar\lambda\bar w}dw\wedge d\bar w\over (w+z)\bar w}=
{i\over 2(2\pi)^2}\int\limits_{w\in\C}
{e^{i(w\bar z+\bar w z)}dw\wedge d\bar w\over w(\bar w-i\lambda)}$$
is called the Faddeev-Green function for the operator
$$\mu\mapsto\bar\pa(\pa+\lambda dz)\mu.$$
From [N2] it follows that\ $\forall\lambda\in\C\b E$ the function
$\psi\big|_{b\Omega}$ can be found through Dirichlet-to-Neumann mapping by
integral equation
$$\psi(z,\lambda)\big|_{b\Omega}=e^{\lambda z}+\int\limits_{\xi\in b\Omega}
e^{\lambda(z-\xi)}g(z-\xi,\lambda)
(\hat\Phi\psi(\xi,\lambda)-\hat\Phi_0\psi(\xi,\lambda)),\eqno(0.4)$$
where  $\hat\Phi\psi=\bar\pa\psi\big|_{b\Omega}$,
$\hat\Phi_0\psi=\bar\pa\psi_0\big|_{b\Omega}$,
$\psi_0\big|_{b\Omega}=\psi\big|_{b\Omega}$ and
$\pa\bar\pa\psi_0\big|_{\Omega}=0$.

By results of [BC1], [GN] and [N2] it follows that $\psi(z,\lambda)$
satisfies $\bar\pa$-equation of Bers-Vekua type with respect to
$\lambda\in\C\b E$:
$$\eqalignno{
&{\pa\psi\over \pa\bar\lambda}=b(\lambda)\bar\psi,\ \ {\rm where}&(0.5)\cr
&\bar\lambda b(\lambda)=-{1\over 2\pi i}\int\limits_{z\in b\Omega}
e^{\lambda z-\bar\lambda\bar z}\bar\pa_z\mu(z,\lambda)=
{1\over 4\pi}\int\limits_{\Omega}e^{\lambda z-\bar\lambda\bar z}q\mu,&(0.6)\cr
&\psi(z,\lambda)e^{-\lambda z}=\mu(z,\lambda)\to 1,\ \lambda\to\infty,\
\forall\ z\in\C.&(0.7)\cr}$$
From  [BC2] and [Na] it follows that for
$q={dd^c\sqrt{\sigma}\over \sqrt{\sigma}}$, $\sigma>0$,
$\sigma\in C^{(2)}(\bar\Omega)$ the exceptional set $E=\{\emptyset\}$ and
function $\lambda\mapsto b(\lambda)$ belongs to
$L^{2+\ep}(\C)\cap L^{2-\ep}(\C)$ for some $\ep>0$.
As a consequence function $\mu=e^{-\lambda z}\psi$ is a unique solution of the
 Fredholm integral equation
$$\mu(z,\lambda)+{1\over 2\pi}\int\limits_{\lambda^{\prime}\in\C}
b(\lambda^{\prime})e^{\bar\lambda^{\prime}\bar z-\lambda^{\prime}z}
\bar\mu(z,\lambda^{\prime})
{d\lambda^{\prime}\wedge d\bar\lambda^{\prime}\over
{\lambda^{\prime}-\lambda}}=1.\eqno(0.8)$$

Integral equations (0.4), (0.8) permit, starting from the
Dirichlet-to-Neumann mapping, to  find firstly the boundary values
$\psi\big|_{b\Omega}$, secondly "$\bar\pa$-scattering data" $b(\lambda)$ and
thirdly function $\psi\big|_{\Omega}$. From equality
$dd^c\psi={dd^c\sqrt{\sigma}\over \sqrt{\sigma}}\psi$ on $X$ we find
finally  ${dd^c\sqrt{\sigma}\over \sqrt{\sigma}}$ on $X$.

The scheme of the Bukhgeim type [B] can be presented in the following way.
Let $q=Qdd^c|z|^2$, where $Q\in C^{(1)}(\bar\Omega)$, but potential $Q$ is
not necessary of the conductivity form
${dd^c\sqrt{\sigma}\over \sqrt{\sigma}}$. By variation of Faddeev
statement and proof we obtain that \ $\forall\ a\in\C\ \exists$ compact set
$E\subset\C$ such that\ $\forall\lambda\in\C\b E$ there exists a unique
solution $\psi_a(z,\lambda)$ of the equation $dd^c\psi=q\psi$ with asymptotics
$$\psi_a(z,\lambda)e^{-\lambda(z-a)^2}=\mu_a(z,\lambda)=1+o(1),\ z\to\infty.$$
Such a solution can be found from integral equation (0.3), where kernel
$g(z,\lambda)$ is replaced by kernel
$$g_a(z,\zeta,\lambda)={ie^{\lambda a^2-\bar\lambda\bar a^2}\over
2\pi^2}\int\limits_{\C}
{e^{-\lambda(\zeta-\eta+a)^2+\bar\lambda(\bar\zeta-\bar\eta+\bar a)^2}\over
(\eta-z)(\bar\zeta-\bar\eta)}d\eta\wedge d\bar\eta.$$
Kernel $g_a(z,\zeta,\lambda)$ can be called the Faddeev type Green function
for the operator

\noindent
$\mu\to\bar\pa(\pa+\lambda d(z-a)^2)\mu$.
Equation $\bar\pa(\pa+\lambda d(z-a)^2)\mu={i\over 2}q\mu$ and Green formula
implies
$$\int\limits_{b\Omega}e^{\lambda(z-a)^2-\bar\lambda(\bar z-\bar a)^2}
\bar\pa\mu=
\int\limits_{\Omega}e^{\lambda(z-a)^2-\bar\lambda(\bar z-\bar a)^2}
{q\mu\over 2i}.\eqno(0.9)$$
Stationary phase method, applied to the integral in the right-hand side of
(0.9), gives for $\tau\to\infty$, $\tau\in\R$,  equality
$$\lim_{\scriptstyle \tau\to\infty \atop \tau\in\R}{4\tau\over \pi i}
\int\limits_{z\in b\Omega}e^{i\tau [(z-a)^2+(\bar z-\bar a)^2]}
\bar\pa_z\mu_a(z,i\tau)=Q(a).\eqno(0.10)$$
Formula (0.10) means that values of potential $Q$ in the arbitrary point $a$
of
$\Omega$ can be reconstructed from Dirichlet-to-Neumann mapping
$\mu_a\big|_{b\Omega}\mapsto \bar\pa_z\mu_a\big|_{b\Omega}$
for family of functions $\mu_a(z,\lambda)$ depending on parameter
$\lambda=i\tau$, $\tau>const$, where we assume that $\mu_a\big|_{b\Omega}$
is found using an analog of (0.4) for $\psi_a\big|_{b\Omega}$.

Bukhgeim's scheme works well at least \ $\forall\ Q\in C^{(1)}(\bar\Omega)$.

More constructive scheme of [N2] works quite well only in the absence of
exceptional set $E$ in the $\lambda$-plane for Faddeev type functions.
In papers [BLMP], [Ts], [N3] it was constructed modified Faddeev-Green
 function permitting to
solve inverse boundary problem (0.2), on the $\R^2=\C$,
at least,  under some smallness
assumptions on potential $Q$.

Let us note that the first uniqueness results in the two-dimensional inverse
boundary value or scattering problems for (0.1) or (0.2) goes back to
A.Calderon [C], V.Druskin [D], R.Kohn, M.Vogelius [KV], J.Sylvester, G.Uhlmann
[SU] and R.Novikov [N1].

Note in this connection that the first seminal results on reconstruction of
the two-dimensional Schr\"odinger operator $H$ on the torus from the
data "extracted" from the family of eigenfunctions (Bloch-Floquet) of  single
 energy level $H\psi=E\psi$ were obtained in series of papers
starting from B.Dubrovin, I.Krichever, S.P.Novikov [DKN], S.P.Novikov,
A.Veselov [NV]. These results were obtained in connection with (2+1)-
dimensional evolution equations.

This paper extends to the case of Riemann surfaces reconstruction procedures
 of [N2] and  of [B]. The paper extends (and also corrects) the
recent paper [HM2] where the inverse boundary value problem on
Riemann surface was firstly considered.Earlier in [HM1] it was proved that
if $X\subset\R^3$ possesses a constant conductivity then
$X$ with complex structure can be effectively reconstructed by at most
three generic potential $\to$ current measurements on $bX$.

Very recently, motivated by [B] and [HM1], [HM2], C.Guillarmou and L.Tzou [GT]
 have obtained general identifiability result (without reconstruction
procedure): if for all solutions of equations $dd^cu+q_ju=0$,
$q_j\in C^{(2)}(X)$, $j=1,2$, Cauchy datas $u\big|_{bX}$, $d^cu\big|_{bX}$,
coincide, then $q_1=q_2$ on $X$.

\vskip 2 mm
{\bf 1. Preliminaries and main results}

Let ${\C}P^3$ be complex projective space with homogeneous coordinates
\vfill\eject

\noindent
$w=(w_0:w_1:w_2:w_3)$.
Let ${\C}P^2_{\infty}=\{w\in {\C}P^3:\ w_0=0\}$.
Then ${\C}P^3\b {\C}P^2_{\infty}$ can be considered as
 the complex affine space with coordinates $z_k=w_k/w_0$, $k=1,2,3$.
By classical result of G. Halphen (see R.Hartshorne [H], ch.IV, $\S$ 6) any
compact Riemann surface of genus $g$ can be embedded in ${\C}P^3$ as
projective algebraic curve $\tilde V$, which intersects ${\C}P^2_{\infty}$
transversally in $d>g$ points, where $d\ge 1$ if $g=0$, $d\ge 3$ if $g=1$
and $d\ge g+3$ if $g\ge 2$. Without loss of generality one can suppose that
\item{  i)} $V=\tilde V\b {\C}P^2_{\infty}$ is connected affine algebraic
curve in $\C^3$ defined by polynomial equations
$V=\{z\in\C^3:\ p_1(z)=p_2(z)=p_3(z)=0\}$ such that the rang of the matrix
$[{\pa p_1\over \pa z}(z),{\pa p_2\over \pa z}(z),{\pa p_3\over \pa z}(z)]
\equiv 2\ \forall\ z\in V$.
\item{ ii)} $\tilde V\cap {\C}P^2_{\infty}=\{\beta_1,\ldots,\beta_d\}$, where
$$\beta_l=(0:\beta_l^1:\beta_l^2:\beta_l^3),\
\bigl({\beta_l^2\over \beta_l^1},{\beta_l^3\over \beta_l^1}\bigr)\in
\C^2,\ l=1,2,\ldots,d.$$
\item{iii)} For $r_0>0$ large enough
$$\det\left|\matrix{
{\pa p_{\alpha}\over \pa z_2}\ &\  {\pa p_{\alpha}\over \pa z_3} \hfill\cr
{\pa p_{\beta}\over \pa z_2}\ &\  {\pa p_{\beta}\over \pa z_3} \hfill\cr}
\right|\ne 0\ \
{\rm for}\ \ z\in V:\ |z_1|\ge r_0\ \ {\rm and}\ \ \alpha\ne\beta.$$
\item{ iv)} For $|z|$ large enough:
$${dz_2\over dz_1}\big|_{V_l}=\gamma_l+{\gamma_l^0\over z_1^2}+
O\bigl({1\over z_1^3}\bigr),\ \
{dz_3\over dz_1}\big|_{V_l}=\tilde\gamma_l+{\tilde\gamma_l^0\over z_1^2}+
O\bigl({1\over z_1^3}\bigr),$$
where $\gamma_l,\ \tilde\gamma_l,\ \gamma_l^0,\ \tilde\gamma_l^0\ne 0$,
for $l=1,\ldots,d$, $d\ge 2$.

Let $V_0=\{z\in V:\ |z_1|\le r_0\}$ and $V\b V_0=\cup_{l=1}^dV_l$, where
$\{V_l\}$ are connected components of $V\b V_0$.
Let us equip $V$ by Euclidean volume form $d d^c|z|^2$. Let

\noindent
$\tilde W^{1,\tilde p}(V)=\{F\in L^{\infty}(V):\
\bar\pa F\in L_{0,1}^{\tilde p}(V)\}$,
$\tilde W^{1,\tilde p}_{1,0}(V)=\{f\in L^{\infty}_{1,0}(V):\
\bar\pa f\in L_{1,1}^{\tilde p}(V)\}$, $\tilde p>2$.
Let $H_{0,1}(V)$ denotes the space of antiholomorphic (0,1)-forms on $V$.
Let

\noindent
$H^p_{0,1}(V)=H_{0,1}(V)\cap L_{0,1}^p(V)$, $1<p<2$.

Let $W^{1,p}(V)=\{F\in L^p(V):\ \bar\pa F\in L^p_{0,1}(V)\}$.

From the Hodge-Riemann decomposition theorem (see [GH], [Ho])\
$\forall\Phi_0\in W^{1,p}_{0,1}(\tilde V)$ we have

\noindent
$\Phi_0=\bar\pa(\bar\pa^*G\Phi_0)+{\cal H}\Phi_0$, where
${\cal H}\Phi_0\in H_{0,1}(\tilde V)$ and $G$ is the Hodge-Green operator
for the Laplacian $\bar\pa\bar\pa^*+\bar\pa^*\bar\pa$ on $\tilde V$ with the
properties: $G(H_{0,1}(\tilde V))=0$, $\bar\pa G=G\bar\pa$,
$\bar\pa^*G=G\bar\pa^*$.

Straight generalization of Proposition 1 from
[He] gives explicit operators:

\noindent
$R_1:\ L_{0,1}^p(V)\to L^{\tilde p}(V)$,
$R_0:\ L_{0,1}^p(V)\to\tilde W^{1,\tilde p}(V)$ and
${\cal H}:\ L_{0,1}^p(V)\to H_{0,1}^p(V)$, $1<p<2$,
${1\over \tilde p}={1\over p}-{1\over 2}$, such that\
$\forall\Phi\in L_{0,1}^p(V)$ we have decomposition of Hodge-Riemann type:
$$\eqalign{
&\Phi=\bar\pa R\Phi+{\cal H}\Phi,\ \ {\rm where}\ \ R=R_1+R_0,\cr
&R_1\Phi(z)={1\over 2\pi i}\int\limits_{\xi\in V}\Phi(\xi)\wedge
(dp_{\alpha}\wedge dp_{\beta})\rfloor d\xi_1\wedge d\xi_2\wedge d\xi_3
\det[{\pa p_{\alpha}(\xi)\over \pa\xi},{\pa p_{\beta}(\xi)\over \pa\xi},
{{\bar\xi-\bar z}\over |\xi-z|^2}],\cr
&R_0\Phi(z)=(\bar\pa^*G(\bar\pa R_1\Phi-\Phi))(z)-
(\bar\pa^*G(\bar\pa R_1\Phi-\Phi))(\beta_1),\cr
&(\bar\pa R_1\Phi-\Phi)\in W^{1,p}_{0,1}(\tilde V),\
G\ \  {\rm is\ the\ Hodge-Green\
operator\ for\ Laplacian}\ \ \bar\pa\bar\pa^*\cr
&{\rm for}\ \ (0,1)-{\rm forms}\ \ {\rm on}\ \ \tilde V,\cr}$$
(1,1)-form under sign of integral does not depend on the choice of indexes
$\alpha,\beta=1,2,3$, $\alpha\ne\beta$,
$${\cal H}\Phi=\sum\limits_{j=1}^g
\bigl(\int\limits_V\Phi\wedge\omega_j\bigr)\bar\omega_j,$$
$\{\omega_j\}$
is orthonormal basis of holomorphic (1,0)-forms on $\tilde V$, i.e.
$$\int\limits_V\omega_j\wedge\bar\omega_k=\delta_{jk},\ \ j,k=1,2,\ldots,g.$$

Note that as a corollary of construction of $R$ we have that
$\lim\limits_{\scriptstyle z\in V_1 \atop z\to\infty}
R\Phi(z)=R\Phi(\beta_1)=0$.

\vskip 2 mm
{\bf Remark 1.1.}
If $V=\{z\in\C^2:\ P(z)=0\}$ be algebraic curve in $\C^2$ then formula for
operator $R_1$ is reduced to the following:
$$R_1\Phi(z)={1\over 2\pi i}\int\limits_{\xi\in V}\Phi(\xi){d\xi_1\over
{\pa P\over \pa\xi_2}}\det\bigl[{\pa P\over \pa\xi}(\xi),
{{\bar\xi-\bar z}\over |\xi-z|^2}\bigr].$$

\vskip 2 mm
{\bf Remark 1.2.}
Based on [HP] one can construct an explicit formula not only for the main
part $R_1$ of the $R$-operator, but for the whole operator $R=R_1+R_0$.

Let  $\v\in L^1_{1,1}(V)\cap L_{1,1}^{\infty}(V)$,
$f\in\tilde W_{1,0}^{1,\tilde p}(V)$, $\lambda\in\C$, $\theta\in\C$.

Let
$$\eqalign{
&\hat R_{\theta}\v=R((dz_1+\theta dz_2)\rfloor\v)(dz_1+\theta dz_2),\cr
&R_{\lambda,\theta}f=
e_{-\lambda,\theta}\overline{R(\overline{e_{\lambda,\theta}f})},\ \ {\rm
where}\ \
e_{\lambda,\theta}(z)=
e^{\lambda(z_1+\theta z_2)-\bar\lambda(\bar z_1+\bar\theta\bar z_2)}.\cr}$$
By straight generalization of Propositions 2, 3 from [He] the form
$f=\hat R_{\theta}\v$ is a solution of $\bar\pa f=\v$ on $V$, function
$u=R_{\lambda,\theta}f$ is a solution of
$$\eqalign{
&(\pa+\lambda(dz_1+\theta dz_2))u=f-{\cal H}_{\lambda,\theta}f,\ \ {\rm where}
\cr
&{\cal H}_{\lambda,\theta}f\buildrel \rm def \over = e_{-\lambda,\theta}
\overline{{\cal H}(\overline{e_{\lambda,\theta}f})},\ \
u\in W^{1,\tilde p}(V),\ \ \tilde p>2.\cr}$$
In addition, by straight generalization of Proposition 4 from [He] we have
that
$$\bar\pa(\pa+\lambda(dz_1+\theta dz_2))u=
\v+\bar\lambda(d\bar z_1+\bar\theta d\bar z_2)
\wedge {\cal H}_{\lambda,\theta}(\hat R_{\theta}\v)\ \ {\rm on}\ \ V.$$

\vskip 2 mm
{\bf Definition 1.1.}
The kernel $g_{\lambda,\theta}(z,\xi)$, $z,\xi\in V$, $\lambda\in\C$, of
integral operator $R_{\lambda,\theta}\circ\hat R_{\theta}$ is called in
[He] the Faddeev type Green function for operator
$\bar\pa(\pa+\lambda(dz_1+\theta dz_2))$.

\vskip 2 mm
{\bf Definition 1.2.}
Let $g=genus\,\tilde V$. Let $\{\omega_j\}$, $j=1,\ldots,g$, be orthonormal
basis of holomorphic forms on $\tilde V$. Let $\{a_1,\ldots,a_g\}$ be
different points (or effective divisor) on $V\b V_0$.
Let
$$\Delta_{\theta}(\lambda)=\det\bigl[\int\limits_{\xi\in V}
\hat R_{\theta}(\delta(\xi,a_j))\wedge\bar\omega_k(\xi)
e_{\lambda,\theta}(\xi),\ \ j,k=1,\ldots,g,\bigr]$$
where $\delta(\xi,a_j)$- Dirac (1,1)-form concentrated in $\{a_j\}$.

Let $E_{\theta}=\{\lambda\in\C:\ \Delta_{\theta}(\lambda)=0\}$.

\vskip 2 mm
{\bf Definition 1.3.}
Parameter $\theta\in\C$ will be called generic if
$\theta\notin \{\theta_1,\ldots,\theta_d\}$, where $\theta_l=-1/\gamma_l$.
Divisor $\{a_1,\ldots,a_g\}$ on $V\b V_0$ will be called generic if
$$\det\bigl[{\omega_j\over dz_1}(a_k)\bigr]_{j,k=1,\ldots,g}
\ne 0.$$

\vskip 2 mm
{\bf Proposition 1.1.}
{\sl
Let parameter $\theta\in\C$ and  divisor $\{a_1,\ldots,a_g\}$ on $V\b V_0$
be generic, where
$V_0=\{z\in V:\ |z_1|\le r_0\}$, $g\ge 1$. Then for $r_0$
large enough we have inequalities:
$$\eqalign{
&\overline{\lim}_{\lambda\to\infty}
|\lambda^g\Delta_{\theta}(\lambda)|<\infty\ \ {\it and}\cr
&\forall\ep>0\ \  \underline{\lim}_{\lambda\to\infty}
|\lambda^g\Delta_{\theta}(\lambda)|_{\ep}>0,\ \ {\it where}\ \
|\lambda^g\Delta_{\theta}(\lambda)|_{\ep}=
\sup\limits_{\{\lambda^{\prime}: |\lambda^{\prime}-\lambda|<\ep\}}
|(\lambda^{\prime})^g\cdot\Delta_{\theta}(\lambda^{\prime})|\cr}$$
Besides, the set $E_{\theta}$
is a closed nowhere dense subset of $\C$.
}

Let $X$ be domain containing $V_0$ and relatively compact on $V$. Let
$\sigma\in C^{(3)}(V)$, $\sigma>0$, on $V$, $\sigma=1$ on $V\b X$. Let
$Y$ be domain containing $\bar X$ and relatively compact on $V$. Let divisor
$\{a_1,\ldots,a_g\}$  on $Y\b X$ and parameter $\theta\in\C$ be  generic.

\vskip 2 mm
{\bf Definition 1.4.}
The functions $\psi_{\theta}(z,\lambda)=\sqrt{\sigma}
F_{\theta}(z,\lambda)=\mu_{\theta}(z,\lambda)e^{\lambda(z_1+\theta z_2)}$,
$z\in V$, $\theta\in\C\b \{\theta_1,\ldots,\theta_d\}$,
$\lambda\in\C\b E_{\theta}$, will be called the
Faddeev type functions, associated with $\sigma$, $\theta$ and
$\{a_1,\ldots,a_g\}$ if $\psi_{\theta}$, $F_{\theta}$, $\mu_{\theta}$
satisfy correspondingly properties:
$$\eqalign{
&d\sigma d^cF_{\theta}=2\sqrt{\sigma}e^{\lambda(z_1+\theta z_2)}
\sum\limits_{j=1}^gC_{j,\theta}(\lambda)\delta(z,a_j),\cr
&dd^c\psi_{\theta}=q\psi_{\theta}+
2e^{\lambda(z_1+\theta z_2)}
\sum\limits_{j=1}^gC_{j,\theta}(\lambda)\delta(z,a_j),\cr
&\bar\pa(\pa+\lambda(dz_1+\theta dz_2))\mu_{\theta}={i\over 2}q\mu_{\theta}+
i\sum\limits_{j=1}^gC_{j,\theta}(\lambda)\delta(z,a_j),\cr}\eqno(1.1)$$
and  the normalization condition
$$\lim\limits_{\scriptstyle z\in V_1 \atop z\to\infty}\mu_{\theta}(z,\lambda)
=1,\eqno(1.2)$$
where $\mu_{\theta}\big|_Y\in L^{\tilde p}(Y)$,
$\mu_{\theta}\big|_{V\b Y}\in L^{\infty}(V\b Y)$,  $\tilde p>2$,
$q={dd^c\sqrt{\sigma}\over \sqrt{\sigma}}$, $\{C_{j,\theta}\}$ are some
functions of $\lambda\in\C\b E_{\theta}$.

\vskip 2 mm
{\bf Theorem 1.1.}
{\sl
Under the aforementioned notations and conditions,\ $\forall$ generic $\theta\in\C$, \
$\forall$ generic divisor
$\{a_1,\ldots,a_g\}\subset V\b X$ and \
$\forall$ $\lambda\in\C\b E_{\theta}:\ |\lambda|>
const(V,\{a_j\},\theta,\sigma)$
there exists unique Faddeev type function
$$\psi_{\theta}(z,\lambda)=\sqrt{\sigma}F_{\theta}(z,\lambda)=
e^{\lambda(z_1+\theta z_2)}\mu_{\theta}(z,\lambda),$$
associated with conductivity function $\sigma$ and  divisor
$\{a_1,\ldots,a_g\}$. Moreover:

A) function $z\to\psi_{\theta}(z,\lambda)$ and parameters
$\{C_{j,\theta}(\lambda)\}$ can be found from the following equations,
depending on parameters $\theta\in\C$, $\lambda\in\C\b E_{\theta}$,
$$\eqalign{
&\psi_{\theta}(z,\lambda)-{i\over 2}\int\limits_{\xi\in X}
e^{\lambda((z_1-\xi_1)+\theta(z_2-\xi_2))}g_{\lambda,\theta}(z,\xi)
{dd^c\sqrt{\sigma}\over \sqrt{\sigma}}\psi_{\theta}(z,\lambda)=\cr
&e^{\lambda(z_1+\theta z_2)}+
i\sum\limits_{j=1}^gC_{j,\theta}(\lambda)g_{\lambda,\theta}(z,a_j)
e^{\lambda(z_1+\theta z_2)},\cr}\eqno(1.3)$$
$$\eqalign{
&2\sum\limits_{j=1}^gC_{j,\theta}(\lambda)e_{\lambda,\theta}(a_j)
{\bar\omega_k\over {d\bar z_1+\bar\theta d\bar z_2}}(a_j)=\cr
&-\int\limits_{z\in V}
e^{-\bar\lambda(\bar z_1+\bar\theta\bar z_2)}
{dd^c\sqrt{\sigma}\over \sqrt{\sigma}}\psi_{\theta}(z,\lambda)
{\bar\omega_k\over {d\bar z_1+\bar\theta d\bar z_2}}(z),\cr}\eqno(1.4)$$
where $k=1,2,\ldots,g$ and $\{\omega_j\}$ is orthonormal basis of
holomorphic forms on $\tilde V$;

B) functions $z\to\psi_{\theta}(z,\lambda)$ and parameters
$\{C_{j,\theta}(\lambda)\}$ satisfy the following properties for
$\lambda\in\C\b E_{\theta}:\ |\lambda|\ge const(V,\{a_j\},\theta,\sigma)$

$$\eqalignno{
&\exists\ \lim_{\scriptstyle z\to\infty,\ z\in V_l \atop l=1,2,\ldots,d}
{{\bar z_1+\bar\theta\bar z_2}\over \bar\lambda}
e^{-\bar\lambda(\bar z_1+\bar\theta\bar z_2)}
\bigl({\pa\psi_{\theta}\over \pa\bar z_1}+
\bar\theta {\pa\psi_{\theta}\over \pa\bar z_2}\bigr)\buildrel \rm \over =
\lim_{\scriptstyle z\to\infty \atop z\in V_l}\psi_{\theta}
e^{-\lambda(z_1+\theta z_2)}b_{\theta}(\lambda),&(1.5)\cr
&iC_{j,\theta}(\lambda)=(2\pi i) Res_{a_j}e^{-\lambda(z_1+\theta z_2)}
\pa\psi_{\theta}\buildrel \rm def \over =2\pi i\lim\limits_{\ep\to 0}
\int\limits_{|z-a_j|=\ep}e^{-\lambda(z_1+\theta z_2)}\pa\psi_{\theta},
&(1.6)\cr
&{\pa\psi_{\theta}(z,\lambda)\over \pa\bar\lambda}=b_{\theta}(\lambda)
\overline{\psi_{\theta}(z,\lambda)},&(1.7)\cr
&{\pa C_{j,\theta}(\lambda)\over \pa\bar\lambda}
e^{\lambda(a_{j,1}+\theta a_{j,2})}=
b_{\theta}(\lambda)\overline{C_{j,\theta}(\lambda)}
e^{\bar\lambda(\bar a_{j,1}+\bar\theta \bar a_{j,2})}.&(1.8)\cr}$$
Besides,
$$\eqalign{
&\bar\lambda b_{\theta}(\lambda)d=
-{1\over 2\pi i}\int\limits_{z\in bX}
e_{\lambda,\theta}(z)\bar\pa\mu(z)+i\sum_{j=1}^gC_{j,\theta}
e_{\lambda,\theta}(a_j),\cr
&|\lambda|\cdot |b_{\theta}(\lambda)|\le const(V,\{a_j\},\sigma)
{1\over (|\lambda|+1)^{1/3}}
{1\over |\Delta_{\theta}(\lambda)|(1+|\lambda|)^g},\cr
&|C_{j,\theta}(\lambda)|\le const(V,\{a_j\},\sigma)
{1\over (|\lambda|+1)^{1/3}}
{1\over |\Delta_{\theta}(\lambda)|(1+|\lambda|)^g}.\cr}\eqno(1.9)$$
}

\vskip 2 mm
{\bf Remark 1.3.}
If $\|\ln\sqrt{\sigma}\|_{C^{(2)}(X)}\le const(V,\{a_j\},\theta)$ then the
condition

\noindent
$\lambda\in\C\b E_{\theta}:\ |\lambda|\ge const(V,\{a_j\},\theta,\sigma)$
in Theorem 1.1 can be replaced by the condition $\lambda\in\C\b E_{\theta}$.
Dependence of $const(V,\{a_j\},\theta,\sigma)$ of $\sigma$ means its
dependence only of $\|\ln\sqrt{\sigma}\|_{C^{(2)}(X)}$.

\vskip 2 mm
{\bf Definition 1.5.}
The functions $b_{\theta}(\lambda)$ and $\{C_{j,\theta}\}$
will be called "scattering" data  for potential $q$.
\vskip 2 mm
Let  $\hat\Phi(\psi\big|_{bX})=\bar\pa\psi\big|_{bX}$ for all sufficiently
regular solutions $\psi$ of (0.2) in $\bar X$, where
$q={dd^c\sqrt{\sigma}\over \sqrt{\sigma}}$. The operator $\Phi$ is 
equivalent to the Dirichlet-to-Neumann operator for (0.1). Let $\hat\Phi_0$
denote $\hat\Phi$ for $q\equiv 0$ on $\bar X$.

\vskip 2 mm
{\bf Theorem 1.2.}
{\sl
Under the conditions of Proposition 1.1  and Theorem 1.1, the following
statements are valid:

A) $\forall\lambda\in\C\b E_{\theta}:\ |\lambda|\ge 
const(V,\{a_j\},\theta,\sigma)$
the restriction of
$\psi_{\theta}(z,\lambda)$ on $bX$ and data $\{C_{j,\theta}(\lambda)\}$
can be reconstructed from Dirichlet-to-Neumann data as
unique  solution of the Fredholm integral equation
$$\eqalignno{
&\psi_{\theta}(z,\lambda)|_{bX}+\int\limits_{\xi\in bX}
e^{\lambda[(z_1-\xi_1)+\theta(z_2-\xi_2)]}g_{\lambda,\theta}(z,\xi)
(\hat\Phi-\hat\Phi_0)\psi_{\theta}(\xi,\lambda)=&(1.10)\cr
&e^{\lambda(z_1+\theta z_2)}+i\sum\limits_{j=1}^gC_{j,\theta}(\lambda)
g_{\lambda,\theta}(z,a_j)e^{\lambda(z_1+\theta z_2)},\ \ {\it where}\cr
&\int\limits_{z\in bX}(z_1+\theta z_2)^{-k}
(\pa+\lambda(dz_1+\theta dz_2))
\mu_{\theta}(z,\lambda)=-\sum\limits_{j=1}^g(a_{j,1}+\theta a_{j,2})^{-k}
C_{j,\theta}(\lambda),&(1.11)\cr}$$
$k=2,\ldots,g+1$, where (without restriction of generality) we suppose that
values $\{a_{j,1}\}$ of the first coordinates of points $\{a_j\}$ are
mutually different;

B) Function $\sigma(w)$, $w\in X$, can be reconstructed from Dirichlet-to-
Neumann data
$$\psi_{\theta}\big|_{bX}\buildrel \rm def \over =
\mu_{\theta}\big|_{bX}e^{\lambda(z_1+\theta z_2)}\to
\bar\pa\psi_{\theta}\big|_{bX}$$
by explicit formulas, where we assume that $\psi_{\theta}\big|_{bX}$ is
found using (1.10), (1.11).

For the case $V=\{z\in\C^2:\ P(z)=0\}$, where $P$ is a polynomial of degree
$N$, this formula has the following form.
Let $\{w_m\}$ be points of $V$, where $(dz_1+\theta dz_2)\big|_V(w_m)=0$,
$m=1,\ldots,M$. Then for almost all $\theta$ values
${d d^c\sqrt{\sigma}\over \sqrt{\sigma}d d^c|z|^2}\big|_V(w_m)$ can be
found from the following linear system
$$\eqalign{
&\tau(1+o(1)){d^k\over d\tau^k}\bigl(\int\limits_{z\in bX}
e_{i\tau,\theta}(z)\bar\pa\mu_{\theta}(z,i\tau)\bigr)=\cr
&\sum_{m=1}^M{i\pi (1+|\theta|^2)\over 2}
{d d^c\sqrt{\sigma}\over \sqrt{\sigma}d d^c|z|^2}\bigg|_V(w_m)\times\cr
&{|{\pa P\over \pa z_1}(w)|^3{d^k\over d\tau^k}\exp i\tau
[(w_{m,1}+\theta w_{m,2})+(\bar w_{m,n}+\bar\theta\bar w_{m,2})]\over
\big|{\pa^2 P\over \pa z_1^2}\bigl({\pa P\over \pa z_2}\bigr)^2-
2{\pa^2 P\over \pa z_1\pa z_2}\bigl({\pa P\over \pa z_2}\bigr)
\bigl({\pa P\over \pa z_1}\bigr)
+{\pa^2 P\over \pa z_2^2}\bigl({\pa P\over \pa z_1}\bigr)^2\big|(w_m)},
\cr}\eqno(1.12)$$
where $m,k=1,\ldots,M$; $M=N(N-1)$, $\tau\in\R$, $\tau\to\infty$,
$|\tau|^g|\Delta_{\theta}(i\tau)|\ge\ep>0$, $\ep$- small enough.
Determinant of system (1.12) is proportional to the determinant of
Vandermonde.

\vskip 2 mm
C) If $g=0$ and if $\theta=\theta(\lambda)=\lambda^{-2}$,
 then \ $\forall\ z\in X$  and \ $\forall\lambda\in\C$ function
$\mu_{\theta}(z,\lambda)=\psi_{\theta}(z,\lambda)e^{-\lambda(z_1+\theta z_2}$
is unique solution of Fredholm integral equation
$$\eqalign{
&\mu_{\theta(\lambda)}(z,\lambda)+{1\over 2\pi i}\int_{\xi\in\C}
b_{\theta(\xi)}(\xi)e^{\bar\xi(\bar z_1+\bar\theta(\xi)\bar z_2)-
\xi(z_1+\theta(\xi) z_2)}
\bar\mu_{\theta(\xi)}(z,\xi){d\xi\wedge d\bar\xi\over {\xi-\lambda}}=1,\cr
&{\it where}\ \ |b_{\theta(\xi)}(\xi)|\le {const(V)\over (1+|\xi|)^2},\cr}$$
and function $z\to\sigma(z)$, $z\in X$, can be found from equality
$$dd^c\psi_{\theta(\lambda)}(z,\lambda)=
{dd^c\sqrt{\sigma}\over \sqrt{\sigma}}(z)
\psi_{\theta(\lambda)}(z,\lambda),\ \ z\in X.$$
}

\vskip 2 mm
{\bf Remark 1.4.}
Using the Faddeev type Green function constructed in [He], in [HM2] were
obtained natural analogues of the main  steps of the reconstruction scheme of
[N2] on the Riemann surface $V$. In particular, under a smallness assumption
on $\pa\log\sqrt{\sigma}$ the existence (and uniqueness) of the solution
$\mu(z,\lambda)$ of the Faddeev type integral equation
$$\mu_{\theta}(z,\lambda)=1+{i\over 2}\int\limits_{\xi\in V}
g_{\lambda,\theta}(z,\xi){\mu_{\theta}(\xi,\lambda)
dd^c\sqrt{\sigma}\over\sqrt{\sigma}}+
i\sum_{j=1}^gC_jg_{\lambda,\theta}(z,a_j),\ z\in V,\ \lambda\in\C,$$
holds for any a priori fixed constants $C_1,\ldots,C_g$.
However (and this fact was overlooked in [HM]) for $\lambda\in\C\b E$
there exists  unique choise of constants $C_j(\lambda,\sigma)$ for
which the integral equation above is equivalent to the differential equation
$$\bar\pa(\pa+\lambda(dz_1+\theta dz_2))\mu-{i\over 2}
\bigl({dd^c\sqrt{\sigma}\over\sqrt{\sigma}}\mu\bigr)+
i\sum_{j=1}^gC_j\delta(z,a_j),$$
where $\delta(z,a_j)$ are  Dirac measures concentrated in the points $a_j$.

\vskip 2 mm
{\bf 2. Faddeev type functions on Riemann surfaces. Uniqueness}

Let projective algebraic curve $\tilde V$ be embedded in ${\C}P^3$ and
intersect

\noindent
${\C}P^2_{\infty}=\{w\in {\C}P^3:\ w_0=0\}$ transversally in
$d>g$ points. Let $V=\tilde V\b {\C}P^2_{\infty}$,

\noindent
$V_0=\{z\in V:\ |z_1|\le r_0\}$ and properties i)-iv) from $\S$ 1 be valid.

\vskip 2 mm
{\bf Proposition 2.1.}
{\sl
Let $\sigma$ be positive function belonging to $C^{(2)}(V)$ such that
$\sigma\equiv const=1$ on
$V\b X\subset V\b V_0=\cup_{l=1}^dV_l$, where
$\{V_l\}$ are connected components of $V\b \bar V_0$. Put
$q={dd^c\sqrt{\sigma}\over \sqrt{\sigma}}$. Let $\{a_1,\ldots,a_g\}$ be
 generic divisor with support in $Y\b\bar X$,
$\bar X\subset Y\subset\bar Y\subset V$. Let for generic
$\theta\in\C$ and $\lambda\in\C:\ |\lambda|\ge const(V,\{a_j\},\theta,\sigma)$
 function
$z\mapsto\mu=\mu_{\theta}(z,\lambda)$ be such that:
$$\eqalign{
&\mu\big|_Y\in L^{\tilde p}(Y),\ \ \mu\big|_{V\b Y}\in L^{\infty}(V\b\bar Y),
\cr
&\bar\pa\mu\big|_Y\in L^p(Y),\ \ \bar\pa\mu\big|_{V\b \bar Y}\in
L^{\tilde p}(V\b Y),\ \ 1\le p<2,\ \tilde p>2,\cr}\eqno(2.1)$$
$$\eqalignno{
&\bar\pa(\pa+\lambda(dz_1+\theta dz_2)\mu={i\over 2}q\mu+i\sum_{j=1}^g
C_j\delta(z,a_j)\ \ {\it with\ some}\ \ C_j=C_{j,\theta}(\lambda)\ \ 
{\it and}&(2.2)\cr
&\mu_{\theta}(z,\lambda)\to 0,\ \ z\to\infty,\ \ z\in V_1.&(2.3)\cr}$$
Then $\mu_{\theta}(z,\lambda)\equiv 0$, $z\in V$.
}

\vskip 2 mm
{\bf Remark 2.1.}
Proposition 2.1 is a corrected version of Proposition 2.1 of [HM2]. For the
case $V=\C$ the equivalent result goes back to [BC2].

\vskip 2 mm
{\bf Lemma 2.1.}
{\sl
Let $\psi=\sqrt{\sigma}F=e^{\lambda(z_1+\theta z_2)}\mu$, where $\mu$
satisfies (2.1), (2.2) and
$$F_1=\sqrt{\sigma}\pa F,\ \ F_2=\sqrt{\sigma}\bar\pa F.\eqno(2.4)$$
Then forms $F_1$, $F_2$ satisfy the system of equations
$$\eqalign{
&\bar\pa F_1+F_2\wedge\pa\ln\sqrt{\sigma}=ie^{\lambda(z_1+\theta z_2)}
\sum_{j=1}^gC_j\delta(z,a_j),\cr
&\pa F_2+F_1\wedge\bar\pa\ln\sqrt{\sigma}=-ie^{\lambda(z_1+\theta z_2)}
\sum_{j=1}^gC_j\delta(z,a_j).\cr}\eqno(2.5)$$
}
\vskip 2 mm
{\it Proof of Lemma 2.1.}
From definition of $F_1$ and $F_2$ it follows that
$$\eqalign{
&d\sigma d^cF=i[2\sigma\pa\bar\pa F-\bar\pa\sigma\wedge\pa F+
\pa\sigma\wedge\bar\pa F]=\cr
&2i\sqrt{\sigma}(\pa F_2+F_1\wedge\bar\pa\ln\sqrt{\sigma})=
-2i\sqrt{\sigma}(\bar\pa F_1+F_2\wedge\pa\ln\sqrt{\sigma}).\cr}$$
From (2.4) and (2.2) we deduce also that
$$d(\sigma d^cF)=\sqrt{\sigma}\bigl(dd^c\psi-
\psi{dd^c\sqrt{\sigma}\over \sqrt{\sigma}}\bigr)=
2\sqrt{\sigma}e^{\lambda(z_1+\theta z_2)}\sum_{j=1}^gC_j\delta(z,a_j).$$
These equalities imply (2.5).

Lemma 2.1 is proved.

\vskip 2 mm
{\bf Lemma 2.2.}
{\sl
Let $\{b_m\}$ be the points of $X$, where $(dz_1+\theta dz_2)\big|_X(b_m)=0$.

Let $B^0=\cup_m\{b_m\}$ and $A^0=\cup_j\{a_j\}$.

Let $u_{\pm}=m_1\pm e_{-\lambda,\theta}(z)\bar m_2$, where
$m_1=e^{-\lambda(z_1+\theta z_2)}f_1$, $m_2=e^{-\lambda(z_1+\theta z_2)}f_2$,

\noindent
$f_1=\sqrt{\sigma}{\pa F\over \pa z_1}$,
$f_2=\sqrt{\sigma}{\pa F\over \pa \bar z_1}$.
Let also $q_1={\pa\ln\sqrt{\sigma}\over \pa z_1}$ and
$\delta_0(z,a_j)={\delta(z,a_j)\over d z_1\wedge d\bar z_1}$. Then in
conditions of Lemma 2.1
$$\eqalign{
&\sup\limits_{z\in X}|\bar\pa u_{\pm}\big|_X(z)\cdot dist^2(z,B^0)|=
O\bigl(\sup\limits_{z\in X}|u_{\pm}dist(z,B^0)|\bigr)<\infty;\cr
&u_{\pm}\big|_{V\b X}\in L^1(V\b X)\cap O(V\b (X\cup A^0))\cr}\eqno(2.6)$$
and system (2.5) is equivalent to the system
$$\eqalign{
&{\pa u_{\pm}\over \pa\bar z_1}d\bar z_1=\mp(e_{-\lambda,\theta}(z)
q_1\bar u_{\pm})d\bar z_1+\cr
&i\sum_{j=1}^g(C_j\pm\bar C_je_{-\lambda,\theta}(z))\delta_0(z,a_j)d\bar z_1.
\cr}\eqno(2.7)$$
}

\vskip 2 mm
{\it Proof of Lemma 2.2.}
From (2.1) we deduce the property
$$\eqalign{
&u_{\pm}\big|_Y\in L^p(Y),\ \ 1\le p<2,\cr
&u_{\pm}\big|_{V\b Y}\in L^{\tilde p}(V\b Y)\oplus L^{\infty}(V\b Y),\ \
\tilde p>2.\cr}$$
System (2.5) is equivalent to the system
of equations
$$\eqalign{
&{\pa f_1\over \pa\bar z_1}=-f_2q_1+ie^{\lambda(z_1+\theta z_2)}
\sum_{j=1}^gC_j\delta_0(z,a_j),\cr
&{\pa f_2\over \pa z_1}=-f_1\bar q_1+ie^{\lambda(z_1+\theta z_2)}
\sum_{j=1}^gC_j\delta_0(z,a_j).\cr}$$
This system and definition of $m_1$, $m_2$ imply
$$\eqalign{
&{\pa m_1\over \pa\bar z_1}=-q_1m_2+i\sum_{j=1}^gC_j\delta_0(z,a_j),\cr
&{\pa m_2\over \pa z_1}+\lambda m_2\bigl(1+\theta{\pa z_2\over \pa z_1}\bigr)
=-\bar q_1m_1+i\sum_{j=1}^gC_j\delta_0(z,a_j).\cr}$$

From the last equalities and definition of $u_{\pm}$ we deduce
$$\eqalign{
&{\pa u_{\pm}\over \pa\bar z_1}={\pa m_1\over \pa\bar z_1}\pm
e_{-\lambda,\theta}(z)\bigl({\pa\bar m_2\over \pa\bar z_1}+\bar\lambda
\bigl(1+\bar\theta{\pa\bar z_2\over \pa\bar z_1}\bigr)\bar m_2\bigr)=
-q_1m_2+i\sum_{j=1}^gC_j\delta_0(z,a_j)\pm\cr
&e_{-\lambda,\theta}(z)\bigl(\bar\lambda
\bigl(1+\bar\theta{\pa\bar z_2\over \pa\bar z_1}\bigr)\bar m_2-
\bar\lambda\bar m_2\bigl(1+\bar\theta{\pa\bar z_2\over \pa\bar z_1}\bigr)-
q_1\bar m_1+i\sum_{j=1}^g\bar C_j\delta_0(z,a_j)\bigr)=\cr
&\mp(e_{-\lambda,\theta}(z)q_1\bar u_{\pm})+
i\sum_{j=1}^g(C_j\pm\bar C_je_{-\lambda,\theta}(z))\delta_0(z,a_j).\cr}$$

Property (2.7) is proved.

For proving (2.6) we will use construction coming back to Bers and Vekua
(see [Ro], [V]).
Let $\beta_{\pm}$ be continuous on $Y$ solutions of $\bar\pa$- equations
$$\bar\pa\beta_{\pm}=\pm e_{-\lambda,\theta}(z)q_1{\bar u_{\pm}\over u_{\pm}}
d\bar z_1,$$
where the right-hand side belongs to $L_{0,1}^{\infty}(Y)$.

Functions $v_{\pm}=u_{\pm}e^{-\beta_{\pm}}$ belongs to ${\cal O}(Y)$.
Indeed, from (2.1), (2.2) it follows that
$\mu\in W^{1,p}(Y)\cap W^{1,\tilde p}_{loc}(Y\b (A^0\cup B^0))$.
From this and from definition of $v_{\pm}$ we deduce that
$\bar\pa v_{\pm}=q_1\bar u_{\pm}d\bar z_1e^{-\beta_{\pm}}-
q_1u_{\pm}{\bar u_{\pm}\over u_{\pm}}e^{-\beta_{\pm}}d\bar z_1=0$ on
$Y\b (A^0\cup B^0)$ and the
following formula for $u_{\pm}$ is valid
$$u_{\pm}(z)=v_{\pm}(z)e^{\beta_{\pm}(z)}.\eqno(2.8)$$
From this and (2.7), (2.8) we obtain (2.6).

Lemma 2.2 is proved.

\vskip 2 mm
{\bf Lemma 2.3.}
{\sl
Let $u_{\pm}$ be the functions from Lemma 2.2 and $\mu$ be the function from
Lemma 2.1. Then
$$u_{\pm}={\pa\mu\over \pa z_1}+\lambda
\bigl(1+\theta{\pa z_2\over \pa z_1}\bigr)\mu-q_1\mu\pm
e_{-\lambda,\theta}(z)\bigl({\pa\bar\mu\over \pa z_1}-q_1\bar\mu\bigr).$$
}

\vskip 2 mm
{\it Proof of Lemma 2.3.}
We have
$$u_{\pm}=e^{-\lambda(z_1+\theta z_2)}f_1\pm e^{-\lambda(z_1+\theta z_2)}
\bar f_2=e^{-\lambda(z_1+\theta z_2)}(f_1\pm \bar f_2),$$
where
$$\eqalign{
&f_1=\sqrt{\sigma}{\pa F\over \pa z_1}=\sqrt{\sigma}{\pa\over \pa z_1}
\bigl({1\over \sqrt{\sigma}}e^{\lambda(z_1+\theta z_2)}\mu\bigr)=\cr
&e^{\lambda(z_1+\theta z_2)}
\bigl({\pa\mu\over \pa z_1}+\lambda\bigl(1+\theta{\pa z_2\over \pa z_1}\bigr)
\mu-q_1\mu\bigr),\cr
&\bar f_2=\sqrt{\sigma}{\pa\bar F\over \pa z_1}=
\sqrt{\sigma}{\pa\over \pa z_1}
\bigl({1\over \sqrt{\sigma}}e^{\bar\lambda(\bar z_1+\bar\theta\bar z_2)}
\bar\mu\bigr)=\cr
&e^{\bar\lambda(\bar z_1+\bar\theta \bar z_2)}
\bigl({\pa\bar\mu\over \pa z_1}-q_1\bar\mu\bigr).\cr}$$

This imply Lemma 2.3.

\vskip 2 mm
{\bf Lemma 2.4.}
{\sl
Let $\omega_1,\ldots,\omega_g$ be orthonormal basis of holomorphic 1-forms on
$\tilde V$.
Let $\{a_1,\ldots,a_g\}$ be generic divisor on $Y\b\bar X$, where
$V_0\subset\bar X\subset Y\subset V$.
Put $\omega^0_{j,k}={\omega_k\over dz_1}(a_j)$. Let for some
generic $\theta\in\C$ and $\lambda\in\C$
functions  $u_{\pm}$ from Lemmas 2.2-2.3 satisfy (2.6), (2.7) with some
$C_j=C_{j,\theta}(\lambda)$. Then
$$\eqalign{
&\sup\limits_j|C_{j,\theta}(\lambda)|\le const(V,\{a_j\},\theta)
\|\ln\sqrt{\sigma}\|^2_{W^{2,\infty}(X)}
(1+|\lambda|)^{-1/3}\|u_{\pm}\|_{L^{\infty}(X,B^0)},\cr
&{\it where}\ \ \|u_{\pm}\|_{L^{\infty}(X,B^0)}\buildrel \rm def \over =
\sup\limits_{z\in X}|u_{\pm}(z)dist(z,B^0)|.\cr}$$
}

\vskip 2 mm
{\bf Proof of Lemma 2.4.}
From condition iv) of section 1 we deduce  $|\omega^0_{j,k}|<\infty$.
From definition of generic divisor we obtain $\det[\omega^0_{j,k}]\ne 0$.
From (2.7) and from definition of Dirac measure\ $\forall\ k=1,\ldots,g$
we deduce
$$\eqalign{
&\overline{\lim_{r\to\infty}}\bigl(\int\limits_{\{z\in V:\ |z_1|=r\}}
u_{\pm}\wedge\omega_k\bigr)\pm
\int\limits_Xe_{-\lambda,\theta}(z){\pa\ln\sqrt{\sigma}\over
\pa z_1}\bar u_{\pm}d\bar z_1\wedge\omega_k=\cr
&i\int\limits_Y\sum_{j=1}^g(C_j\pm\bar C_je_{-\lambda,\theta}(z))
\delta_0(z,a_j)d\bar z_1\wedge\omega_k=\cr
&i\sum_{j=1}^g(C_j\pm\bar C_je_{-\lambda,\theta}(a_j))\omega^0_{j,k},\
j,k=1,2,\ldots, g.\cr}\eqno(2.9)$$

From estimates
$\overline{\lim\limits_{r_n\to\infty}}
\sup\limits_{\{z\in V:\ |z_1|=r_n\}}|u_{\pm}(z)|<\infty$, for some sequence
$r_n\to\infty$, and
${|\omega_k|\over dz_1}\le O\bigl(|{1\over z_1^2}|\bigr)$, $z\in V\b Y$,
$k=1,\ldots,g$, we obtain
$$\overline{\lim_{r\to\infty}}
\bigg|\int\limits_{\{z\in V:\ |z_1|=r\}}u_{\pm}\wedge\omega_k\bigg|=0.
\eqno(2.10)$$
From (2.9), (2.10) and Kramers's formula we obtain
$$\eqalign{
&i(C_j\pm\bar C_je_{-\lambda,\theta}(a_j))=\cr
&{\det[\omega^0_{1,k};\ldots;\omega^0_{j-1,k};\int\limits_X\pm
e_{-\lambda,\theta}(z){\pa\ln\sqrt{\sigma}\over \pa z_1}\bar u_{\pm}
d\bar z_1\wedge\omega_k;\omega^0_{j+1,k};\ldots;\omega^0_{g,k}]\over
\det[\omega^0_{j,k}]},\cr}\eqno(2.11)$$
where $j,k=1,\ldots,g$.

Let us prove estimate
$$\eqalign{
&\bigg|\int\limits_X
e_{-\lambda,\theta}(z){\pa\ln\sqrt{\sigma}\over \pa z_1}\bar u_{\pm}
d\bar z_1\wedge\omega_k\bigg|\le\cr
&const(X,\theta)(1+|\lambda|)^{-1/3}
\|\ln\sqrt{\sigma}\|^2_{W^{2,\infty}(X)}\cdot\|u_{\pm}\|_{L^{\infty}(X,B^0)}.
\cr}\eqno(2.12)$$
For $|\lambda|\le 1$ estimate follows directly, using  that
$\ln\sqrt{\sigma}\in W^{1,\infty}(X)$.

Let $B^{\ep}=\cup_{m=1}^M\{z\in X:\ |z-b_m|\le\ep\}$.

Let $\chi_{\ep,\nu}$, $\nu=1,2$, be functions from $C^{(1)}(V)$ such that
$\chi_{\ep,1}+\chi_{\ep,2}\equiv 1$ on $V$,

\noindent
$supp\,\chi_{\ep,1}\subset B^{2\ep}$, $supp\,\chi_{\ep,2}\subset V\b B^{\ep}$,
$|d\chi_{\ep,\nu}|=O({1\over \ep})$, $\nu=1,2$.

Put $J_{\nu}^{\ep}u_{\pm}=\int\limits_X\chi_{\ep,\nu}(z)
e_{-\lambda,\theta}(z){\pa\ln\sqrt{\sigma}\over \pa z_1}\bar u_{\pm}
d\bar z_1\wedge\omega_k$, $\nu=1,2$.
We have directly:
$$|J_1^{\ep}u_{\pm}|\le const(X)\ep
\|\ln\sqrt{\sigma}\|_{W^{1,1}(X)}\cdot \|u_{\pm}\|_{L^{\infty}(X,B^0)}.
\eqno(2.13)$$
For $J_2^{\ep}u_{\pm}$ we obtain by integration by parts:
$$\eqalign{
&J_2^{\ep}u_{\pm}=-{1\over \lambda}\int\limits_X\chi_{\ep,2}
\pa e_{-\lambda,\theta}(z){\pa\ln\sqrt{\sigma}\over \pa z_1}\bar u_{\pm}
d\bar z_1\wedge {\omega_k\over {dz_1+\theta dz_2}}=\cr
&{1\over \lambda}\int\limits_Xe_{-\lambda,\theta}(z)
\pa\biggl(\chi_{\ep,2}{\pa\ln\sqrt{\sigma}\over \pa z_1}\bar u_{\pm}
d\bar z_1\wedge {\omega_k\over {dz_1+\theta dz_2}}\biggr).\cr}\eqno(2.14)$$
To estimate (2.14) we use (2.6) and the following properties:
$|\pa\chi_{\ep,2}|=O({1\over \ep})$,

\noindent
$supp(\pa\chi_{\ep,2})\subset B^{2\ep}$,
$$\eqalign{
&\|{\pa\ln\sqrt{\sigma}\over \pa z_1}d\bar z_1\wedge\pa\chi_{\ep,2}u_{\pm}
{\omega_k\over {dz_1+\theta dz_2}}\|_{L^1_{0,1}(X)}\le\cr
&{const(X,\theta)\over \ep}
\|\ln\sqrt{\sigma}\|_{W^{1,\infty}(X)}\|u_{\pm}\|_{L^{\infty}(X,B^0)}\cr
&\|{\pa^2\ln\sqrt{\sigma}\over \pa z_1^2}dz_1\wedge d\bar z_1\chi_{\ep,2}
u_{\pm}{\omega_k\over {dz_1+\theta dz_2}}\|_{L^1_{1,1}(X)}\le\cr
&|\ln\ep|const(X,\theta)\|\ln\sqrt{\sigma}\|_{W^{2,\infty}(X)}
\|u_{\pm}\|_{L^{\infty}(X,B^0)}\cr
&\|{\pa\ln\sqrt{\sigma}\over \pa z_1}d\bar z_1\chi_{\ep,2}u_{\pm}\wedge
\pa\bigl({\omega_k\over {dz_1+\theta dz_2}}\bigr)\|_{L^1_{0,1}(X)}\le\cr
&{const(X,\theta)\over \ep}
\|\ln\sqrt{\sigma}\|_{W^{1,\infty}(X)}\|u_{\pm}\|_{L^{\infty}(X,B^0)}\cr
&\pa\bar u_{\pm}\big|_X=\mp(e_{\lambda,\theta}(z)\bar q_1\bar u_{\pm})dz_1.
\cr}$$
From (2.14), (2.6)  and these properties we obtain
$$\eqalign{
&|J_2^{\ep}u_{\pm}|\le |\ln\ep|
{const(X,\theta)\over |\lambda|}\|\ln\sqrt{\sigma}\|_{W^{2,\infty}(X)}\cdot
\|u_{\pm}\|_{L^{\infty}(X,B^0)}+\cr
&{const(X,\theta)\over \ep |\lambda|}\|\ln\sqrt{\sigma}\|_{W^{1,\infty}(X)}
\cdot\|u_{\pm}\|_{L^{\infty}(X,B^0)}+\cr
&{const(X,\theta,\delta)\over \ep^{1+\delta}|\lambda|}
\|\ln\sqrt{\sigma}\|_{W^{1,\infty}(X)}\cdot
\|u_{\pm}\|_{L^{\infty}(X,B^0)}.\cr}\eqno(2.15)$$
Putting in (2.13), (2.15) $\ep={1\over \sqrt{\lambda}}$ and $\delta=1/3$
we obtain (2.12) for
$|\lambda|\ge 1$.

Inequalities (2.11), (2.12) imply estimate
$$|C_j\pm\bar C_je_{-\lambda,\theta}(a_j)|\le const(X,\{a_j\},\theta)
(1+|\lambda|)^{-1/3}\|\ln\sqrt{\sigma}\|^2_{W^{2,\infty}(X)}
\cdot\|u_{\pm}\|_{L^{\infty}(X,B^0)}.$$
We obtained statement of Lemma 2.4.

\vskip 2 mm
{\bf Lemma 2.5.}
{\sl
Let functions $u_{\pm}$ satisfy (2.6), (2.7) and $R$ - operator from section
1. Then
$$\|R[e_{-\lambda,\theta}q_1\bar u_{\pm}d\bar\xi_1\|_{L^{\infty}(X,B^0)}\le
const(X,\theta)(1+|\lambda|)^{-1/5}
\|\ln\sqrt{\sigma}\|_{W^{2,\infty}(X)}\cdot
\|u_{\pm}\|_{L^{\infty}(X,B^0)}.$$
}

\vskip 2 mm
{\it Proof of Lemma 2.5.}

Let $\chi_{\ep,\nu}$, $\nu=1,2$,  be partition of unity from Lemma 2.4.
Put $S_{\nu}^{\ep}u_{\pm}=R[\chi_{\ep,\nu}q_1\bar u_{\pm}d\bar\xi_1]$,
$\nu=1,2$. Using (2.6) and formula for operator $R$ we deduce estimate
$$\|S_1^{\ep}u_{\pm}\|_{L^{\infty}(X,B^0)}=
O(\ep)\|\ln\sqrt{\sigma}\|_{W^{1,\infty}(X)}
\|u_{\pm}\|_{L^{\infty}(X,B^0)}.\eqno(2.16)$$
Let $R_{1,0}(\xi,z)$ be kernel of operator $R$. It means, in particular, that
 $\bar\pa_{\xi}R_{1,0}(\xi,z)=-\delta(\xi,z)$, where $\delta(\xi,z)$- Dirac
(1,1)- measure, concentrated in the point $\xi=z$.
We have
$$S_2^{\ep}u_{\pm}=\int\limits_X\chi_{\ep,2}e_{-\lambda,\theta}
q_1\bar u_{\pm}d\bar\xi_1R_{1,0}(\xi,z).\eqno(2.17)$$
Integration by parts in (2.17) gives the following
$$\eqalign{
&S_2^{\ep}u_{\pm}={1\over \bar\lambda}
\int\limits_X\bar\pa e_{-\lambda,\theta}(\xi){d\bar\xi_1\over
{d\bar\xi_1+\bar\theta d\bar\xi_2}}
\chi_{\ep,2}(\xi)
q_1(\xi)\bar u_{\pm}(\xi)R_{1,0}(\xi,z)=\cr
&-{1\over \bar\lambda}
\int\limits_Xe_{-\lambda,\theta}(\xi)\bar\pa\biggl({d\bar\xi_1\over
{d\bar\xi_1+\bar\theta d\bar\xi_2}}\chi_{\ep,2}(\xi)
q_1(\xi)\bar u_{\pm}(\xi)\biggr)R_{1,0}(\xi,z)+\cr
&{1\over \bar\lambda}
e_{-\lambda,\theta}(z){d\bar\xi_1\over
{d\bar\xi_1+\bar\theta d\bar\xi_2}}(z)\chi_{\ep,2}(z)
q_1(z)\bar u_{\pm}(z).\cr}\eqno(2.18)$$
To estimate (2.18) we use (2.6), properties  of partition of unity
$\{\chi_{\ep,\nu}\}$ and inequalities
$$\eqalign{
&\bigg|{d\bar\xi_1\over {d\bar\xi_1+\bar\theta d\bar\xi_2}}(\xi)\bigg|=
O\bigl({1\over dist(\xi,B^0)}\bigr),\ \
\bigg|\bar\pa{d\bar\xi_1\over {d\bar\xi_1+\bar\theta d\bar\xi_2}}(\xi)\bigg|=
O\bigl({1\over (dist(\xi,B^0))^2}\bigr),\cr
&|q_1(\xi)|=O\bigl({1\over dist(\xi,B^0)}\bigr),\ \
|\bar\pa q_1(\xi)|=O\bigl({1\over (dist(\xi,B^0))^2}\bigr), \ \
\xi\in X.\cr}\eqno(2.19)$$
From (2.19), (2.8) and from the formula for operator $R$ we deduce estimate
$$\|S_2^{\ep}u_{\pm}\|_{L^{\infty}(X)}=O({1\over \ep^4|\lambda|})
\|\ln\sqrt{\sigma}\|_{W^{2,\infty}(X)}
\|u_{\pm}\|_{L^{\infty}(X,B^0)}.\eqno(2.20)$$
Putting in (2.16), (2.20) $\ep={1\over |\lambda|^{1/5}}$ we obtain
statement of Lemma 2.5.

\vskip 2 mm
{\it Proof of Proposition 2.1.}

Let function $\mu$ satisfy conditions (2.1)-(2.3) and $u_{\pm}$ be functions
defined in

\noindent
Lemma 2.2. Then by Lemma 2.3 we have
$$\eqalign{
&\lim_{\scriptstyle z\to\infty \atop z\in V_1}u_{\pm}(z,\lambda)=
\lim_{\scriptstyle z\to\infty \atop z\in V_1}
(m_1\pm e_{-\lambda,\theta}(z)\bar m_2)=\cr
&\lim_{\scriptstyle z\to\infty \atop z\in V_1}
[\lambda\bigl(1+\theta{dz_2\over dz_1}\bigr)\mu+{\pa\mu\over \pa z_1}\pm
e_{-\lambda,\theta}(z){\pa\bar\mu\over \pa z_1}]\to 0.\cr}\eqno(2.21)$$
Let
$$h_{\pm}=u_{\pm}\pm R[(e_{-\lambda,\theta}(z)q_1\bar u_{\pm})d\bar z_1-
i\sum_{j=1}^g(C_j\pm\bar C_je_{-\lambda,\theta}(z))\delta_0(z,a_j)d\bar z_1],
\eqno(2.22)$$
where $R$ is the operator from section 1.

By Lemmas 2.2-2.5 and properties of operator $R$ we have
$h_{\pm}\in {\cal O}(V)\cap L^{\infty}(V)$ and

\noindent
$h_{\pm}(z,\lambda)\to 0$,
$z\to\infty$, $z\in V_1$. By Liouville theorem, $h_{\pm}(z,\lambda)\equiv 0$
on $V$, $\lambda\in\C$. Then from (2.22) with $h_{\pm}(z,\lambda)\equiv 0$
and Lemmas 2.4, 2.5 it follows that $u_{\pm}(z,\lambda)\equiv 0$,
$z\in V$, if $\lambda\in\C\b E_{\theta}:\ |\lambda|\ge const(V,\{a_j\},\theta)
\|\ln\sqrt{\sigma}\|^2_{W^{2,\infty}(X)}$.
Property $u_{\pm}(z,\lambda)\equiv 0$, $z\in V$, implies by Lemma 2.3 equality
 ${\pa\mu\over \pa\bar z_1}-\bar q_1\mu=0$, $z\in V$, where $\mu(z)\to\infty$
if $z\in V_1$, $z\to\infty$.
The Liouville type theorem for generalized holomorphic functions
([Ro], theorem 7.1) implies $\mu\equiv 0$. Proposition 2.1 is proved.

\vskip 2 mm
{\bf 3. Faddeev type functions on Riemann surface. Existence.}

\noindent
{\bf Proof of Theorem 1.1A}

\vskip 2 mm
{\bf Proposition 3.1.}
{\sl
Let conductivity  $\sigma$ and divisor $\{a_1,\ldots,a_g\}$ satisfy
conditions of Proposition 2.1. Then \
$\forall$ generic $\theta\in\C$ and\
$\forall\lambda\in\C\b E_{\theta}:\ |\lambda|\ge
const (V,\{a_j\},\theta,\sigma)$
there exists unique Faddeev type function
$$\eqalign{
&\psi\buildrel \rm def \over =\sqrt{\sigma}F\buildrel \rm def \over =
e^{\lambda(z_1+\theta z_2)}\mu,\ \ {\rm where}\cr
&\psi=\psi_{\theta}(z,\lambda),\ F=F_{\theta}(z,\lambda),\
\mu=\mu_{\theta}(z,\lambda),\cr}\eqno(3.1)$$
associated with $\sigma$ and divisor $\{a_1,\ldots,a_g\}$, i.e.
$$\eqalign{
&\bar\pa(\pa+\lambda(dz_1+\theta dz_2))\mu={i\over 2}q\mu+
\sum_{j=1}^gC_j\delta(z,a_j),\ \ {\it for\ some}\ \
C_j=C_{j,\theta}(\lambda),\ \
{\it where}\cr
&q={dd^c\sqrt{\sigma}\over \sqrt{\sigma}},\ \mu\big|_Y\in L^{\tilde p}(Y),\
\mu\big|_{V\b \bar Y}\in L^{\infty}(V\b \bar Y),\
\lim_{\scriptstyle z\to\infty \atop z\in V_1}\mu_{\theta}(z,\lambda)=1.\cr}
\eqno(3.1a)$$
In addition,
$$\eqalign{
&\|\mu_{\theta}(z,\lambda)-\mu_{\theta}(\infty_l,\lambda)\|_{L^{\tilde p}(V)}
\le {const(V,\{a_j\},\theta,\sigma,\tilde p,\ep)\over
|\Delta_{\theta}(\lambda)|\cdot(1+|\lambda|)^{g+1-\ep}},\cr
&{\it where}\ \ \ \mu_{\theta}(\infty_l,\lambda)
\buildrel \rm def \over =\lim\limits_{\scriptstyle z\to\infty \atop z\in V_l}
\mu_{\theta}(z,\lambda),\ \ l=1,\ldots,d,\cr
&\|\pa\mu\|_{L^p_{1,0}(Y)}+
\|\pa\mu\|_{L^{\tilde p}_{1,0}(V\b Y)}\le
{const(V,\{a_j\},\theta,\sigma,p,\tilde p,\ep)\over
|\Delta_{\theta}(\lambda)|\cdot
(1+|\lambda|)^{g-\ep}},\ \ p<2,\ \tilde p>2,\cr}\eqno(3.1b)$$
$$\eqalign{
&\forall\ \ {\it generic}\ \ \theta\in\C\ \ {\it and}\ \
\lambda\in\C\b E_{\theta}:\ |\lambda|\ge const(V,\{a_j\},\theta,\sigma),\cr
&{\pa\mu\over \pa\bar\lambda}\big|_Y\in W^{1,p}(Y),\ \
{\pa\mu\over \pa\bar\lambda}\big|_{V_l\b Y}
\in L^{\infty}(V_l\b Y)\cup W^{1,\tilde p}(V_l\b Y),\cr}\eqno(3.1c)$$
where $\{V_l\}$ are connected components of $V\b V_0$, $l=1,\ldots,d$,
$$e_{\lambda,\theta}(z)=e^{\lambda(z_1+\theta z_2)-
\bar\lambda(\bar z_1+\bar\theta \bar z_2)}.$$
}

\vskip 2 mm
{\bf Remark 3.1.}
Proposition 3.1 is a corrected version of Proposition 2.2 from [HM2].
For the case $V=\C$  the results of such a type goes back to [F1], [F2].

\vskip 2 mm
{\bf Lemma 3.1.}
{\sl
Under the conditions of Proposition 3.1,\
$\forall\lambda\in\C\b E_{\theta}$ function

\noindent
$z\to\mu_{\theta}(z,\lambda)$
belonging to $L^{\tilde p}(Y)$ on $Y$ and to $L^{\infty}(V\b Y)$ on
$V\b Y$ satisfies (3.1a) iff

there exists $C_j=C_{j,\theta}(\lambda)$,
$j=1,\ldots,g$, such that
$$\mu_{\theta}(z,\lambda)=1+{i\over 2}\int\limits_{\xi\in X}
g_{\lambda,\theta}(z,\xi)q\mu_{\theta}(\xi,\lambda)+
i\sum_{j=1}^gC_{j,\theta}(\lambda)g_{\lambda,\theta}(z,a_j)\eqno(3.2)$$
and one of two equivalent conditions is valid
$$\eqalign{
&{\cal H}_{\lambda,\theta}\bigl(\hat R_{\theta}
\bigl({i\over 2}q\mu\bigr)\bigr)+
i\sum_{j=1}^gC_{j,\theta}(\lambda)
{\cal H}_{\lambda,\theta}(\hat R_{\theta}(\delta(z,a_j))=0\ \ {\rm or}\cr
&(\pa+\lambda(dz_1+\theta dz_2))\mu_{\theta}(z,\lambda)\in
H_{1,0}(V\b (X\cup_{j=1}^g\{a_j\}))\cap L^1_{1,0}(Y\b X),\cr}\eqno(3.3)$$
where $g_{\lambda,\theta}$ is  Faddeev type Green function, $\hat R_{\theta}$,
 ${\cal H}_{\lambda,\theta}$ - operators defined in section 1.
}

\vskip 2 mm
{\it Proof of Lemma 3.1.}
From Proposition 4 in [He] and from definition of Green function
$g_{\lambda,\theta}(z,\xi)$ we deduce that integral equation (3.2) is
equivalent to the following differential equation
$$\eqalign{
&\bar\pa(\pa+\lambda(dz_1+\theta dz_2))\mu={i\over 2}q\mu+
i\sum_{j=1}^gC_{j,\theta}\delta(z,a_j)+\cr
&\bar\lambda(d\bar z_1+\bar\theta d\bar z_2))\times
\big[{\cal H}_{\lambda,\theta}\bigl(\hat R_{\theta}
\bigl({i\over 2}q\mu\bigr)\bigr)+
i\sum_{j=1}^gC_{j,\theta}
{\cal H}_{\lambda,\theta}(\hat R_{\theta}(\delta(z,a_j))\big].\cr}\eqno(3.4)$$
Equation (3.4) is equivalent to (3.1a) if one of two equivalent
conditions (3.3) is valid.

Lemma 3.1 is proved.

\vskip 2 mm
{\bf Lemma 3.2.}
{\sl
Let $\{a_1,\ldots,a_g\}$ be generic divisor in $Y\b \bar X$. Then for any\
 generic $\theta\in\C$ and

\noindent
$\forall\lambda\in\C\b E_{\theta}:\ |\lambda|\ge
const(V,\{a_j\},\theta,\sigma)$,
integral equation
(3.2), (3.3) is uniquely solvable Fredholm integral equation in the space
$\tilde W^{1,\tilde p}(V)$.
}

\vskip 2 mm
{\it Proof of Lemma 3.2.}
Let $\theta\in\C$ and
$\lambda\in\C\b E_{\theta}:\ |\lambda|\ge const(V,\{a_j\},\theta,\sigma)$.
From (3.2), (3.3) we
obtain integral equation for $\tilde\mu_{\theta}=\mu_{\theta}-1$ and
$\tilde C_{j,\theta}$:
$$\eqalign{
&\tilde\mu_{\theta}(z,\lambda)-{i\over 2}\int\limits_{\xi\in V}
g_{\lambda,\theta}(z,\xi)q(\xi)\tilde\mu_{\theta}(\xi,\lambda)-
i\sum_{j=1}^g\tilde C_{j,\theta}(\lambda)g_{\lambda,\theta}(z,a_j)=\cr
&{i\over 2}\int\limits_{\xi\in V}
g_{\lambda,\theta}(z,\xi)q(\xi)+
i\sum_{j=1}^gC_{j,\theta}^0(\lambda)g_{\lambda,\theta}(z,a_j).\cr}\eqno(3.5)$$
Parameters $\tilde C_j=\tilde C_{j,\theta}(\lambda)$, $j=1,\ldots,g$,
are defined by the equations:
$$\eqalign{
&-i\sum_{j=1}^g\tilde C_j\int\limits_V\hat R_{\theta}(\delta(\xi,a_j))
\bar\omega_k(\xi)e_{\lambda,\theta}(\xi)=\cr
&\int\limits_{\xi\in V}e_{\lambda,\theta}(\xi)\hat R_{\theta}
\bigl({i\over 2}q\tilde\mu\bigr)\bar\omega_k(\xi),\ \ k=1,2,\ldots,g.\cr}
\eqno(3.6)$$
We remind that determinant of system (3.6) is exactly
$\Delta_{\theta}(\lambda)$.

Parameters $C_{j,\theta}^0$ are defined by (3.6) with
$C_{j,\theta}^0$ in place of  $\tilde C_{j,\theta}$ and 1
 in place of $\tilde\mu$. One can see also that
$C_{j,\theta}^0(\lambda)=C_{j,\theta}(\lambda)-\tilde C_{j,\theta}(\lambda)$.

Let us prove that (3.5), (3.6) determine Fredholm integral equation in the
space $\tilde W^{1,\tilde p}(V)$, $\tilde p>2$.

Propositions 2, 3 of [He] imply that correspondance
$$\tilde\mu\mapsto R_{\lambda,\theta}\circ(\hat R_{\theta}
\bigl({i\over 2}q\tilde\mu\bigr)+
i\sum_{j=1}^g\tilde C_{j,\theta}\hat R_{\theta}(\delta(z,a_j))\bigr)$$
define linear continuous mapping of $\tilde W^{1,\tilde p}(V)$ into
itself. This mapping is compact because mapping $\tilde\mu\to q\tilde\mu$,
$supp\,q\subset X$, from $\tilde W^{1,\tilde p}(V)$ into
$L^{\tilde p}_{1,1}(X)$ is compact, operator

\noindent
$\hat R_{\theta}:\ L^{\tilde p}_{1,1}(X)\to\tilde W^{1,\tilde p}_{1,0}(V)$ and
operator
$R_{\lambda,\theta}:\ \tilde W^{1,\tilde p}_{1,0}(V)\to
\tilde W^{1,\tilde p}(V)$   are bounded.

If for fixed $\lambda\not\in E_{\theta}$ Fredholm equation
(3.5), (3.6) is not solvable then corresponding homogeneous equation, when
the right-hand side of (3.5) is replaced by zero, admits nontrivial
solution $\tilde\mu^*=\mu^*-1$.

By Lemma 3.1 function  $\tilde\mu^*$ satisfies differential equation (2.2)
with $C_j$ replaced by $\tilde C_j$ and with property $\tilde\mu^*(z)\to 0$,
$z\to\infty$, $z\in V_1$.

By Proposition 2.1, $\tilde\mu^*\equiv 0$ if
$\lambda\in\C\b E_{\theta}:\ |\lambda|\ge const(V,\{a_j\},\theta,\sigma)$.

It means that equation (3.2), (3.3) is uniquely solvable Fredholm integral
equation for any
$\lambda\in\C\b E_{\theta}:\ |\lambda|\ge const(V,\{a_j\},\theta,\sigma)$.

Lemma 3.2 is proved.

\vskip 2 mm
{\bf Lemma 3.3.}
{\sl
Let $\{a_1,\ldots,a_g\}$ be generic divisor on $Y\b X$. Let
$\lambda\in\C\b E_{\theta}$. Let $\mu$ be solution of integral equation
(3.2), (3.3). Then relations (3.3) determining parameters
$C_j=C_{j,\theta}(\lambda)$ are reduced to the following explicit formulas
$$2i\sum_{j=1}^gC_{j,\theta}e_{\lambda,\theta}(a_j)
{\bar\omega_k\over d\bar z_1}(a_j)=
\int\limits_{z\in X}e_{\lambda,\theta}(z)
\bigl(i{dd^c\sqrt{\sigma}\over \sqrt{\sigma}}+2\bar\pa\ln\sqrt{\sigma}\wedge
\pa\ln\sqrt{\sigma}\bigr)\mu{\bar\omega_k\over d\bar z_1}(z).\eqno(3.7)$$
}

\vskip 2 mm
{\it Proof of Lemma 3.3.}
By Lemma 3.1 equations (3.2), (3.3) are equivalent to the equation:
$$\bar\pa(\pa+\lambda(dz_1+\theta dz_2))\mu={i\over 2}q\mu+
i\sum_{j=1}^gC_{j,\theta}\delta(z,a_j),\eqno(3.8)$$
where $\mu=\mu_{\theta}(z,\lambda)\to 1$, $z\in V_1$, $z\to\infty$.

System (2.7) implies the following relation
$$\eqalign{
&\overline{\lim\limits_{R\to\infty}}\int\limits_{|z_1|=R}
\bar u_{\pm}\wedge\bar\omega_k+
i\int\limits_{z\in V\b X}\sum_{j=1}^g(\bar C_{j,\theta}\mp C_{j,\theta}
e_{\lambda,\theta}(z)){\delta(z,a_j)\over d\bar z_1}\bar\omega_k=\cr
&\mp\int\limits_{z\in X}
e^{\lambda(z_1+\theta z_2)-\bar\lambda(\bar z_1+\bar\theta\bar z_2)}
\bar q_1u_{\pm}dz_1\wedge\bar\omega_k,\cr}\eqno(3.9)$$
where $\bar q_1={\pa\ln\sqrt{\sigma}\over \pa\bar z_1}$.

To obtain (3.9) we  multiply the both sides of (2.7) by $\wedge\omega_k$,
integrate on $V$ and  take conjugation.

From Lemma 2.3 and Lemma 3.2 it follows that
$$\eqalign{
&u_{\pm}(z)\to\lambda(1+\theta\gamma_l)\cdot
\lim_{\scriptstyle z\to\infty \atop z\in V_l}\mu_{\theta}(z,\lambda),\
z\to\infty,\ z\in V_l,\cr
&{\rm where}\ \ \gamma_l=
\lim_{\scriptstyle z\to\infty \atop z\in V_l}{\pa z_2\over \pa z_1},\ \
\lim_{\scriptstyle z\to\infty \atop z\in V_1}\mu_{\theta}(z,\lambda)=1.\cr}$$
Existence of
$\lim\limits_{\scriptstyle z\to\infty \atop z\in V_l}\mu_{\theta}(z,\lambda)$
follows from Lemma 4.1 below.
This imply that
$$\overline{\lim\limits_{R\to\infty}}\big|\int\limits_{|z_1|=R}\bar u_{\pm}
\wedge\bar\omega_k\big|=
\overline{\lim\limits_{R\to\infty}}\big|\int\limits_{|z_1|=R}\bar\lambda
(1+\bar\theta\bar\gamma_l)\bar\omega_k\big|=
\lim\limits_{R\to\infty}|\lambda|O\bigl({1\over R}\bigr)=0.\eqno(3.10)$$

From (3.9), (3.10) and definition of $u_{\pm}$ we obtain
$$\eqalign{
&2i\sum_{j=1}^g\int\limits_{V\b X}C_je_{\lambda,\theta}(z)
{\delta(z,a_j)\over d\bar z_1}\wedge\bar\omega_k=
\int\limits_{z\in X}e_{\lambda,\theta}(z)
\bar q_1(u_++u_-)dz_1\wedge\bar\omega_k=\cr
&2\int\limits_{z\in X}e^{-\bar\lambda(\bar z_1+\bar\theta\bar z_2)}
\bar q_1f_1dz_1\wedge\bar\omega_k,\ \ {\rm where}\ \
f_1=\sqrt{\sigma}{\pa F\over \pa z_1}.\cr}$$
By Lemma 2.3 we have
$$\eqalign{
&2\int\limits_{z\in X}e^{-\bar\lambda(\bar z_1+\bar\theta\bar z_2)}
\bar q_1f_1dz_1\wedge\bar\omega_k=\cr
&2\int\limits_{z\in X}e_{\lambda,\theta}(z)\bar q_1
\bigl({\pa\mu\over \pa z_1}+\lambda\mu+\lambda\theta{\pa z_2\over \pa z_1}\mu-
q_1\mu\bigr)dz_1\wedge\bar\omega_k.\cr}\eqno(3.11)$$
From definition of $\delta(z,a_j)$ we have
$$2i\sum_{j=1}^g\int\limits_{z\in V\b X}C_je_{\lambda,\theta}(z)
{\delta(z,a_j)\over d\bar z_1}\wedge\bar\omega_k=
-2i\sum_{j=1}^gC_je_{\lambda,\theta}(a_j)
{\bar\omega_k\over d\bar z_1}(a_j).\eqno(3.12)$$
By integration by part we have
$$\eqalign{
&2\int\limits_{z\in X}e_{\lambda,\theta}(z)\bar q_1
\bigl({\pa\mu\over \pa z_1}+\lambda\mu\bigr)dz_1\wedge\bar\omega_k=
2\int\limits_Xe_{\lambda,\theta}(z)
{\pa\ln\sqrt{\sigma}\over \pa \bar z_1}\lambda\mu dz_1\wedge\bar\omega_k-\cr
&-2\int\limits_Xe_{\lambda,\theta}(z)
{\pa\ln\sqrt{\sigma}\over \pa \bar z_1}\bigl(\lambda\mu+\lambda\theta
{\pa z_2\over \pa z_1}\mu\bigr)dz_1\wedge\bar\omega_k-
2\int\limits_Xe_{\lambda,\theta}(z)
{\pa^2\ln\sqrt{\sigma}\over \pa z_1\pa \bar z_1}\mu dz_1\wedge\bar\omega_k=\cr
&-2\int\limits_Xe_{\lambda,\theta}(z)
\bigl({\pa^2\ln\sqrt{\sigma}\over \pa z_1\pa \bar z_1}+
{\pa\ln\sqrt{\sigma}\over \pa \bar z_1}\lambda\theta
{\pa z_2\over \pa z_1}\bigr)\mu dz_1\wedge\bar\omega_k.\cr}\eqno(3.13)$$
Using (3.11), (3.12), (3.13) we obtain
$$i\sum_{j=1}^gC_{j,\theta}e_{\lambda,\theta}(a_j)
{\bar\omega_k\over d\bar z_1}(a_j)=
\int\limits_{z\in X}e_{\lambda,\theta}(z)
\bigl({\pa^2\ln\sqrt{\sigma}\over \pa z_1\pa \bar z_1}+
\big|{\pa\ln\sqrt{\sigma}\over \pa \bar z_1}\big|^2\bigr)
\mu dz_1\wedge\bar\omega_k.$$

Lemma 3.3 is proved.

\vskip 2 mm
{\it Proof of Proposition 3.1.}
a) By Lemmas 3.1-3.3 statement (3.1a) of Proposition is valid, i.e.
there exists function
$z\to\mu_{\theta}(z,\lambda)$, $z\in V$ with property (3.1a)\
$\forall\lambda\in \C\b E_{\theta}:\ |\lambda|\ge
const(V,\{a_j\},\theta,\sigma)$.

b) Put $f_0=\hat R_{\theta}\bigl({i\over 2}q\mu\bigr)$,
$f_1=\hat R_{\theta}\bigl(i\sum\limits_{j=1}^gC_{j,\theta}
\delta(z,a_j)\bigr)$ and $f=f_0+f_1$.
By (3.2) we have $\mu-1=R_{\lambda,\theta}f=R_{\lambda,\theta}f_0+
R_{\lambda,\theta}f_1$.

Put
$$L^{p,\tilde p}_{0,q}(V)=\{u:\ u\big|_Y\in L^p_{0,q}(Y),\ \
u\big|_{V\b\bar Y}\in L^{\tilde p}_{0,q}(V\b Y)\},\ \ 1\le p<2,\ \tilde p>2,\
q=0,1.$$
By Proposition 3 ii$^{\prime}$ from [He] we obtain
$$\eqalign{
&\|\mu-\mu_{\theta}(\infty_l,\lambda)\|_{L^{\tilde p}(V_l\b Y)}\le\cr
&const(V,\tilde p,\theta)\cdot\min(|\lambda|^{-1/2},|\lambda|^{-1})
\bigl(\|f_0\|_{\tilde W^{1,\tilde p}_{1,0}(V)}+\sum_{j=1}^g|C_{j,\theta}|
\bigr)\cr
&\|\pa\mu\|_{L^{p,\tilde p}_{1,0}(V)}\le const(V,\tilde p,\theta)
\bigl(\|f_0\|_{\tilde W^{1,\tilde p}_{1,0}(V)}+\sum_{j=1}^g|C_{j,\theta}|
\bigr).\cr}\eqno(3.14)$$
For proving estimates (3.1b) let us now estimate $\{C_{j,\theta}^0\}$.

In order to estimate  $\{C_{j,\theta}^0\}$ we must use equations (3.6), where
parameters $\{\tilde C_{j,\theta}\}$ are replaced by  $\{C_{j,\theta}^0\}$
and function
$\tilde\mu$ is replaced by 1. For modified equations (3.6)

1) we  apply Kramer formula for solution of linear system and integration by
parts in all integrals of this system, using
$e_{\lambda,\theta}(z)(d\bar z_1+\bar\theta d\bar z_2)={1\over \bar\lambda}
\bar\pa e_{\lambda,\theta}(z)$. In addition, we use: formula (1.2) for
$\Delta_{\theta}(\lambda)$, formula
$\bar\pa\hat R_{\theta}\bigl({i\over 2}q\mu\bigr)={i\over 2}q\mu$ and
estimate of singular integral, containing
$\bar\pa\bigl({\bar\omega_k\over {d\bar z_1+\bar\theta d\bar z_2}}\bigr)$.
This gives inequality:
$$\sum_j|C_{j,\theta}^0(\lambda)|\le {const(V,\{a_j\},\theta,\sigma)\over
|\Delta_{\theta}(\lambda)|(1+|\lambda|)^g}.$$

ii) The equation (3.5) together with obtained inequality for
$\sum\,|C_{j,\theta}^0(\lambda)|$, estimate of Faddeev type Green function
$|g_{\lambda,\theta}(z,\xi)|=O\bigl({1\over |\lambda|^{1-\ep}}\bigr)$ are
used to obtain estimate (3.15) for $\sum\,|\tilde C_{j,\theta}(\lambda)|$
and $|\mu_{\theta}(\lambda)|$:
$$\eqalign{
&|\lambda|^{-\ep}\|\mu\|_{\tilde W^{1,\tilde p}(V)}+
\sum_j|\tilde C_{j,\theta}(\lambda)|\le
{const(V,\{a_j\},\theta,\sigma,\tilde p,\ep)\over
|\Delta_{\theta}(\lambda)|(1+|\lambda|)^g}\ \ {\rm and}\cr
&\|\mu-\mu(\infty_l,\cdot)\|_{L^{\tilde p}(V_l\b Y)}\le
{const(V,\{a_j\},\theta,\sigma,\tilde p,\ep)\over
|\Delta_{\theta}(\lambda)|(1+|\lambda|)^{g+1-\ep}},\cr
&{\rm where}\ \
\lambda\in \C\b E_{\theta}:\ |\lambda|\ge
const(V,\{a_j\},\theta,\sigma,\ep),\ \ l=1,\ldots,d,\
\mu_{\theta}(\infty_1,\lambda)=1.\cr}\eqno(3.15)$$

These estimates
imply estimates (3.1b).

c) Differentiation of equation (3.2) with respect to $\bar\lambda$ gives
equality
$$\eqalign{
&{\pa\mu\over \pa\bar\lambda}-R_{\lambda,\theta}\circ
\bigl(\hat R_{\theta}\bigl({i\over 2}q{\pa\mu\over \pa\bar\lambda}+
i\sum_{j=1}^g
{\pa C_{j,\theta}(\lambda)\over \pa\bar\lambda}\delta(z,a_j)\bigr)=\cr
&(\bar z_1+\bar\theta\bar z_2)(\mu-1)-R_{\lambda,\theta}
\bigl((\bar\xi_1+\bar\theta\bar\xi_2)\hat R_{\theta}
\bigl({i\over 2}q\mu+i\sum_{j=1}^g
C_{j,\theta}\delta(z,a_j)\bigr).\cr}\eqno(3.16)$$
Equality (3.16) can be rewritten in the following form
$$\eqalign{
&{\pa\mu\over \pa\bar\lambda}=\bigl(I-R_{\lambda,\theta}\circ\hat R_{\theta}
\bigl({i\over 2}q\cdot\bigr)\bigr)^{-1}
\bigl[(\bar z_1+\bar\theta\bar z_2)(\mu-1)+R_{\lambda,\theta}\circ
\hat R_{\theta}\bigl(i\sum_{j=1}^g
{\pa C_{j,\theta}\over \pa\bar\lambda}\delta(z,a_j)\bigr)-\cr
&R_{\lambda,\theta}\bigl((\bar\xi_1+\bar\theta\bar\xi_2)\hat R_{\theta}
\bigl({i\over 2}q\mu+
i\sum_{j=1}^gC_{j,\theta}(\lambda)\delta(z,a_j)\bigr)\bigr].\cr}\eqno(3.17)$$
Using Propositions 2, 3 from [He], estimates from part (b) of this proof
 we obtain from (3.17)
$$\eqalign{
&e_{\lambda,\theta}(z){\pa\mu\over \pa\bar\lambda}\big|_Y\in W^{1,p}(Y),
\cr
&e_{\lambda,\theta}(z){\pa\mu\over \pa\bar\lambda}\big|_{V_l}
\in W^{1,\tilde p}(V_l\b Y)\cup L^{\infty}(V_l\b Y).\cr}$$

Statement (3.1c) is proved.

Proposition 3.1 is proved.

\vskip 2 mm
{\bf 4. Equation ${\pa\mu(z,\lambda)\over \pa\bar\lambda}=b_{\theta}(\lambda)
e_{-\lambda,\theta}(z)\bar\mu_{\theta}(z,\lambda)$.
Proof of Theorem 1.1B}

\vskip 2 mm
{\bf Proposition 4.1.}
{\sl
Let conductivity  $\sigma$, divisor $\{a_1,\ldots,a_g\}$ and $\theta$
satisfy the conditions of Proposition 2.1. Let function
$\psi_{\theta}(z,\lambda)=e^{\lambda(z_1+\theta z_2)}\mu_{\theta}(z,\lambda)$
be the Faddeev type function, associated with $\sigma$, $\theta$ and divisor
$\{a_1,\ldots,a_g\}$. Then for
$\lambda\in \C\b E_{\theta}:\ |\lambda|\ge const(V,\{a_j\},\theta,\sigma)$
\item{ i)} the following $\bar\pa$-equations  take place
$$\eqalignno{
&{\pa\mu_{\theta}(z,\lambda)\over \pa\bar\lambda}=b_{\theta}(\lambda)
e_{-\lambda,\theta}(z)\overline{\mu_{\theta}(z,\lambda)},\ \ {\rm if}\ \
z\in V\b \{a_1,\ldots,a_g\},&(4.1)\cr
&{\pa C_{j,\theta}(\lambda)\over \pa\bar\lambda}=b_{\theta}(\lambda)
e_{-\lambda,\theta}(a_j)\overline{C_{j,\theta}(\lambda)},\ \ j=1,\ldots,g,\ \
 {\rm where}\ \ &(4.2)\cr}$$
\item{ii)} function $b_{\theta}(\lambda)$ satisfies equations:
$$\eqalign{
&b_{\theta}(\lambda)
\lim_{\scriptstyle z\to\infty \atop z\in V_l}
\overline{\mu_{\theta}(z,\lambda)}=
\lim_{\scriptstyle z\to\infty \atop z\in V_l}
{{\bar z_1+\bar\theta\bar z_2}\over \bar\lambda}e_{\lambda,\theta}(z)
{\pa\mu_{\theta}(z,\lambda)\over \pa(\bar z_1+\bar\theta\bar z_2)},\cr
&\bar\lambda b_{\theta}(\lambda)d=-{1\over 2\pi i}\int\limits_{z\in bX}
e_{\lambda,\theta}(z)\bar\pa\mu_{\theta}(z,\lambda)+i\sum_{j=1}^g
C_{j,\theta}(\lambda)e_{\lambda,\theta}(a_j),\ l=1,\ldots,d\cr}\eqno(4.3)$$
and the inequality
$$|\lambda|(1+|\lambda|)^g|\Delta_{\theta}(\lambda)|\cdot|b_{\theta}(\lambda)|\le
const(V,\{a_j\},\theta,\sigma){1\over (|\lambda|+1)^{1/3}}.
\eqno(4.4)$$
}

\vskip 2 mm
{\bf Remark 4.1.}
For the case $V=\C$ this statement is obtained in [GN], [N2], [N3].
Proposition 4.1 is a corrected version of Proposition 3.2 of [HM2].

\vskip 2 mm
{\bf Lemma 4.1.}
{\sl
i) Let function $\mu=\mu_{\theta}(z,\lambda),\ z\in V\b Y$,

\noindent
$\lambda\in\C\b E_{\theta}:\ |\lambda|\ge const(V,\{a_j\},\theta,\sigma)$
satisfy equation
$$\bar\pa(\pa+\lambda(d z_1+\theta dz_2))\mu=0\ \ {\rm on}\ \ V\b Y \eqno(4.5)
$$ and the property
$$\eqalign{
&[\mu-\mu_{\theta}(\infty_l,\lambda)]\big|_{V_l\b Y}\in W^{1,\tilde p}
(V_l\b\bar Y),\ \ {\it where}\ \ \tilde p>2,\cr
&\mu_{\theta}(\infty_l,\lambda)\buildrel \rm def \over =
\lim_{\scriptstyle z\to\infty \atop z\in V_l}
\mu_{\theta}(z,\lambda),\ l=1,\ldots,d.\cr}$$
Then
$$\eqalign{
&A\buildrel \rm def \over ={\pa\mu\over \pa(z_1+\theta z_2)}+
\lambda\mu\in {\cal O}(\tilde V\b\bar Y)\ \ {\it and}\cr
&A\big|_{V_l\b Y}=\lambda\mu(\infty_l)+\sum_{k=1}^{\infty}A_{k,l}
{1\over (z_1+\theta z_2)^k},\cr}$$
$$\eqalign{
&\bar B\buildrel \rm def \over =e_{\lambda,\theta}(z)
{\pa\mu\over \pa(\bar z_1+\bar\theta\bar z_2)}\in
\overline{{\cal O}(\tilde V\b\bar Y)}\ \ {\it and}\cr
&\bar B\big|_{V_l\b Y}=\sum_{k=1}^{\infty}B_{k,l}
{1\over (\bar z_1+\bar\theta\bar z_2)^k},\ l=1,\ldots,d,\cr}\eqno(4.6)$$
where ${\cal O}(\tilde V\b\bar Y)$ is the space of holomorphic functions on
$(\tilde V\b\bar Y)$.

ii) Let
$$M\big|_{V_l}=\mu_{\theta}(\infty_l,\lambda)+\sum\limits_{k=1}^{\infty}
{a_{k,l}(\lambda)\over (z_1+\theta z_2)^k}\ \ {\rm and}\ \
\bar N\big|_{V_l}=\sum\limits_{k=1}^{\infty}
{b_{k,l}(\lambda)\over (\bar z_1+\bar\theta\bar z_2)^k}$$
 be formal series with coefficients determined by relations
$$\lambda a_{k,l}-(k-1)a_{k-1,l}=A_{k,l},\ \
\bar\lambda b_{k,l}-(k-1)b_{k-1,l}=B_{k,l},\ \ l=1,\ldots,d,\ \ k=1,2,\ldots.
$$
Let
$$M_{\nu}\big|_{V_l}=\mu_{\theta}(\infty_l,\lambda)+
\sum_{k=1}^{\nu}{a_{k,l}\over (z_1+\theta z_2)^k},\ \
\bar N_{\nu}\big|_{V_l}=
\sum_{k=1}^{\nu}{b_{k,l}\over (\bar z_1+\bar\theta\bar z_2)^k}.\eqno(4.7)$$
Then function $\mu$ has the asymptotic decomposition
$$\eqalign{
&\mu\big|_{V_l}=M\big|_{V_l}+e_{-\lambda,\theta}(z)\bar N\big|_{V_l},\ \
z_1\to\infty,\cr
&{\rm i.e.}\ \ \mu\big|_{V_l}=M\big|_{V_l}+e_{-\lambda,\theta}(z)
\bar N_{\nu}\big|_{V_l}+O\bigl({1\over |z_1|^{\nu+1}}\bigr).\cr}$$
}

\vskip 2 mm
{\it Proof of Lemma 4.1.}

i) From (4.5) it follows that
$$\pa\bar\pa(e^{\lambda(z_1+\theta z_2)}\mu(z,\lambda))\big|_{V\b\bar Y}=0.$$
Thus
$\bar\pa(e^{\lambda(z_1+\theta z_2)}\mu(z,\lambda))=
e^{\lambda(z_1+\theta z_2)}\bar\pa\mu$ is antiholomorphic form on
$V\b\bar Y$ and $\pa\mu+\lambda\mu(dz_1+\theta dz_2)$ is holomorphic form on
$V\b\bar Y$. From this, condition
$\bar\pa\mu\in L^{\tilde p}_{0,1}(V\b\bar Y)$ and the Cauchy theorem it
follows that
$$\eqalign{
&e^{\lambda(z_1+\theta z_2)}\bar\pa\mu\big|_{V_l\b\bar Y}=
e^{\bar\lambda(\bar z_1+\bar\theta\bar z_2)}
\bar B (d\bar z_1+\bar\theta d\bar z_2)\big|_{V_l\b\bar Y}=\cr
&e^{\bar\lambda(\bar z_1+\bar\theta\bar z_2)}
\sum_{k=1}^{\infty}{B_{k,l}\over (\bar z_1+\bar\theta\bar z_2)^k}
(d\bar z_1+\bar\theta d\bar z_2)\big|_{V_l}\ \ {\rm and}\cr
&(\pa\mu+\lambda\mu(dz_1+\theta dz_2))\big|_{V_l\b\bar Y}=
A(dz_1+\theta dz_2)\big|_{V_l\b\bar Y}=\cr
&\bigl(\lambda\mu(\infty_l)+
\sum_{k=1}^{\infty}{A_{k,l}\over (z_1+\theta z_2)^k}\bigr)
(dz_1+\theta dz_2)\big|_{V_l\b\bar Y}.\cr}$$
It gives (4.6).

ii) From (4.6), (4.7) we obtain, first, that
$$\eqalign{
&\bar\pa\mu\big|_{V_l}=e^{-\lambda(z_1+\theta z_2)}\bar\pa
\bigl(e^{\bar\lambda(\bar z_1+\bar\theta\bar z_2)}\bar N_{\nu}\bigr)
\big|_{V_l}+O\bigl({1\over |\bar z_1|^{\nu+1}}\bigr)\cr
&{\rm then}\ \ \mu\big|_{V_l}=M_{\nu}\big|_{V_l}+e_{-\lambda,\theta}(z)
\bar N_{\nu}\big|_{V_l}+
\tilde O\bigl({1\over |\bar z_1|^{\nu}}\bigr).\cr}\eqno(4.8)$$
Comparison of the last equality  for different indexes $\nu$ and $\nu+1$
implies  that $\tilde O\bigl({1\over |\bar z_1|^{\nu}}\bigr)=
O\bigl({1\over |\bar z_1|^{\nu+1}}\bigr)$.

It gives statement of Lemma 4.1.

\vskip 2 mm
{\bf Lemma 4.2.}
{\sl
i) Functions $M_{\nu}$ and $N_{\nu}$ (congugated to $\bar N_{\nu}$)
from decomposition (4.8) have the
following properties:
$$\eqalign{
&\forall\ z\in\tilde V\b Y\ \exists\ \lim\limits_{\nu\to\infty}
\bigl({\pa M_{\nu}\over \pa(z_1+\theta z_2)}+\lambda M_{\nu}\bigr)
\buildrel \rm def \over ={\pa M\over \pa(z_1+\theta z_2)}+\lambda M\ \
{\it and}\cr
&\exists\ \lim\limits_{\nu\to\infty}
\bigl({\pa N_{\nu}\over \pa(z_1+\theta z_2)}+\lambda N_{\nu}\bigr)
\buildrel \rm def \over ={\pa N\over \pa(z_1+\theta z_2)}+\lambda N.\cr}$$

ii) Functions ${\pa M\over \pa(z_1+\theta z_2)}+\lambda M$ and
${\pa N\over \pa(z_1+\theta z_2)}+\lambda N$ belongs to
${\cal O}(\tilde V\b Y)$ and
$$\eqalign{
&{\pa\mu\over \pa(\bar z_1+\bar\theta\bar z_2)}=e_{-\lambda,\theta}(z)
\bigl({\pa\bar N\over \pa(\bar z_1+\bar\theta\bar z_2)}+
\bar\lambda\bar N\bigr),\cr
&{\pa\mu\over \pa(z_1+\theta z_2)}+\lambda\mu=
{\pa M\over \pa(z_1+\theta z_2)}+\lambda M,\cr}\eqno(4.9)$$
$${\pa N\over \pa(z_1+\theta z_2)}+\lambda N\to 0,\ \ z_1\to\infty.\eqno(4.10)
$$
}

\vskip 2 mm
{\it Proof of Lemma 4.2.}

Part i) and equalities (4.9), (4.10) from part ii) follow directly from (4.8).

Properties (4.8), (4.9), (4.10), property
$\bar\pa\mu\in L^{p,\tilde p}_{0,1}$ (Proposition 3.1b)
and  extension property of bounded holomorphic functions through isolated
singularities imply that
$${\pa M\over \pa(z_1+\theta z_2)}+\lambda M\ \ {\rm and}\ \
{\pa N\over \pa(z_1+\theta z_2)}+\lambda N$$
belongs to ${\cal O}(\tilde V\b Y)$.

Lemma 4.2 is proved.

\vskip 2 mm
{\bf Lemma 4.3.}
{\sl
Let $\psi_{\theta}(z,\lambda)=e^{\lambda(z_1+\theta z_2)}
\mu_{\theta}(z,\lambda)$ be the Faddeev type function on $V$, associated with
potential $q={dd^c\sqrt{\sigma}\over \sqrt{\sigma}}$ and divisor
$\{a_1,\ldots,a_g\}$ on $Y\b\bar X$. Then

\noindent
$\forall\lambda\in \C\b E_{\theta}:\ |\lambda|\ge
const(V,\{a_j\},\theta,\sigma)$
$$\eqalign{
&e_{\lambda,\theta}(z)
{\pa\mu\over \pa(\bar z_1+\bar\theta\bar z_2)}\big|_{V_l\b\bar Y}=
\sum_{k=1}^{\infty}B_{k,l}(\bar z_1+\bar\theta\bar z_2)^{-k},\ \ {\it where}
\cr
&B_{1,l}=-{1\over 2\pi i}\times\cr
&\int\limits_{\{z\in V_l:\ |z_1|=r_1\}}
e_{\lambda,\theta}(z)
{\pa\mu\over \pa(\bar z_1+\bar\theta\bar z_2)}(d\bar z_1+\bar\theta d\bar z_2)
\ \forall\ r_1:\
Y\subset\{z\in V:\ |z_1|<r_1\}.\cr}\eqno(4.11)$$
}

\vskip 2 mm
{\it Proof of Lemma 4.3.}

Estimate of $\pa\mu$ from (3.1b) and the Cauchy theorem, applied
to antiholomorphic function
$e_{\lambda,\theta}(z)
{\pa\mu\over \pa(\bar z_1+\bar\theta\bar z_2)}\big|_{V_l\b\bar Y}$ implies
(4.11).

\vskip 2 mm
{\it Proof of Proposition 4.1.}

Since $\psi$, $\mu$ are Faddeev type functions, we have the equations
$$\eqalign{
&\bar\pa(\pa+\lambda(d z_1+\theta dz_2))\mu={i\over 2}q\mu+
i\sum_{j=1}^{\infty}C_{j,\theta}(\lambda)\delta(z,a_j),\cr
&dd^c\psi=q\psi+2\sum_{j=1}^ge^{\lambda(z_1+\theta z_2)}
C_{j,\theta}(\lambda)\delta(z,a_j).\cr}$$
Put $\psi_{\bar\lambda}={\pa\psi\over \pa\bar\lambda}$ and
$\mu_{\bar\lambda}={\pa\mu\over \pa\bar\lambda}$.

We obtain
$$dd^c\psi_{\bar\lambda}=q\psi_{\bar\lambda}+2\sum_{j=1}^g
e^{\lambda(z_1+\theta z_2)}{\pa C_{j,\theta}\over \pa\bar\lambda}(\lambda)
\delta(z,a_j).$$

From Lemma 4.1 we deduce
$$\eqalign{
&{\pa\mu\over \pa(\bar z_1+\bar\theta\bar z_2)}\big|_{V_l\b\bar Y}
=e_{-\lambda,\theta}(z)
{B_{1,l}(\lambda)\over {\bar z_1+\bar\theta\bar z_2}}+
O\bigl({1\over |z_1|^2}\bigr),\ {\rm and}\cr
&\bigl({\pa\mu\over \pa(z_1+\theta z_2)}+\lambda\mu\bigr)\big|_{V_l\b\bar Y}=
\lambda\mu(\infty_l)+
{A_{1,l}(\lambda)\over {z_1+\theta z_2}}+
O\bigl({1\over |z_1|^2}\bigr).\cr}\eqno(4.12)$$
From (4.6), (4.7), (4.8)  we deduce
$$\eqalignno{
&\mu\big|_{V_l\b\bar Y}=\mu_{\theta}(\infty_l,\lambda)+
{a_l(\lambda)\over {z_1+\theta z_2}}+
e_{-\lambda,\theta}(z){b_l(\lambda)\over {\bar z_1+\bar\theta\bar z_2}}+
O\bigl({1\over |z_1|^2}\bigr),\ z_1\to\infty,&(4.13)\cr
&{\rm where}\ \ \bar\lambda b_l(\lambda)\buildrel \rm def \over =
\bar\lambda b_{1,l}(\lambda)=B_{1,l},\ \
\lambda a_l(\lambda)\buildrel \rm def \over =\lambda a_{1,l}(\lambda)=
A_{1,l},\ \ l=1,\ldots,d.&(4.14)\cr}$$

From (4.13) and (3.1c) we obtain for $l=1,\ldots,d$
$$\eqalign{
&\psi\big|_{V_l\b Y}=e^{\lambda(z_1+\theta z_2)}\mu=\cr
&e^{\lambda(z_1+\theta z_2)}
\bigl(\mu_{\theta}(\infty_l,\lambda)+{a_l(\lambda)\over {z_1+\theta z_2}}+
e^{\bar\lambda(\bar z_1+\bar\theta \bar z_2)-\lambda(z_1+\theta z_2)}
{b_l(\lambda)\over {\bar z_1+\bar\theta\bar z_2}}+
O\bigl({1\over |z_1|^2}\bigr)\bigr),\cr
&\psi_{\bar\lambda}\big|_{V_l\b Y}=
{\pa\psi\over \pa\bar\lambda}\big|_{V_l\b Y}=
e^{\bar\lambda(\bar z_1+\bar\theta \bar z_2)}
\bigl[(\bar z_1+\bar\theta \bar z_2)
{{b_l(\lambda)+e_{\lambda,\theta}(z){\pa\mu_{\theta}(\infty_l,\lambda)\over
\pa\bar\lambda}}
\over {\bar z_1+\bar\theta\bar z_2}}+
O\bigl({1\over |z_1|}\bigr)\bigr]=\cr
&e^{\bar\lambda(\bar z_1+\bar\theta \bar z_2)}
\bigl(b_l(\lambda)+e_{\lambda,\theta}(z){\pa\mu_{\theta}(\infty_l,\lambda)
\over \pa\bar\lambda}+
O\bigl({1\over |z_1|}\bigr)\bigr).\cr}$$

For function $\mu_{\bar\lambda}=
e^{-\lambda(z_1+\theta z_2)}\psi_{\bar\lambda}$ we obtain
$$\eqalign{
&\bar\pa(\pa+\lambda(dz_1+\theta dz_2))\mu_{\bar\lambda}=
{i\over 2}q\mu_{\bar\lambda}+
i\sum_{j=1}^g{\pa C_{j,\theta}\over \pa\bar\lambda}\delta(z,a_j)\cr
&{\rm and}\ \ \mu_{\bar\lambda}=e_{-\lambda,\theta}(z)
\bigl(b_l(\lambda)+e_{\lambda,\theta}(z){\pa\mu_{\theta}(\infty_l,\lambda)
\over \pa\bar\lambda}+
O\bigl({1\over |z_1|}\bigr)\bigr),\ \ z\in V_l.\cr}$$

For $z_1$ large enough function $e_{-\lambda,\theta}(z)\bar\mu_{\bar\lambda}
\buildrel \rm def \over =\v(z,\lambda)$ satisfies equation

\noindent
$\bar\pa(\pa+\lambda(dz_1+\theta dz_2))\v=0$.
From this, Lemma 4.1 and property
$\overline{\lim\limits_{z\to\infty}}|\v(z,\lambda)|_V<\infty$ we deduce
that $\v\big|_{V_l}(z,\lambda)\to const_l(\lambda)\buildrel \rm def \over =
\v(\infty_l,\lambda)$, if $z\in V_l$, $z\to\infty$, $l=1,\ldots,d$.
So in the relations above we have
$e_{\lambda,\theta}(z)\mu_{\bar\lambda}(\infty_l,\lambda)\equiv 0$,
$l=1,\ldots,d$.
Functions  $e_{-\lambda,\theta}(z)\bar\mu_{\bar\lambda}$ and $\mu$ both
satisfy equation $\bar\pa(\pa+\lambda(dz_1+\theta dz_2))\mu={i\over 2}q\mu$
on $V\b \{a_1,\ldots,a_g\}$.
Besides $\overline{\mu\big|_{V_l}(z,\lambda)}\to
\overline{\mu(\infty_l,\lambda)}$
and $e_{\lambda,\theta}(z)\mu_{\bar\lambda}(z,\lambda)\to b_l(\lambda)$, if
$z\in V_l$, $z\to\infty$.
Applying Proposition 2.1 we obtain
$$e_{\lambda,\theta}(z)\mu_{\bar\lambda}=
b_l(\lambda)\bar\mu_{\theta}(z,\lambda)
(\overline{\mu_{\theta}(\infty_l,\lambda)})^{-1},\ \ l=1,\ldots,d.$$

This implies equalities (4.1), (4.2), where
$$b_{\theta}(\lambda)={b_l(\lambda)\over
 \overline{\mu_{\theta}(\infty_l,\lambda)}},\ \ l=1,\ldots,d.
\eqno(4.15)$$
Asymptotic formula (4.3) follows from (4.11), (4.14) and (4.15). These
formulas and Cauchy-Green formula imply also the following important
expression for $b_{\theta}(\lambda)$:
$$\bar\lambda b_{\theta}(\lambda)d=-{1\over 2\pi i}\int\limits_{z\in bY}
e_{\lambda,\theta}(z)\bar\pa\mu=
-{1\over 2\pi i}\int\limits_{z\in bX}
e_{\lambda,\theta}(z)\bar\pa\mu+i\sum_{j=1}^gC_{j,\theta}
e_{\lambda,\theta}(a_j),\eqno(4.16)$$
where
$$\int\limits_{z\in bX}e_{\lambda,\theta}(z)\bar\pa\mu=
\int\limits_X{1\over 2i}e_{\lambda,\theta}(z)q\mu.\eqno(4.17)$$
Equality (4.3) follows from (4.16). This equality together with estimate
of $\{C_j\}$ from Lemma 2.4 and estimate through integration by parts of
$\int\limits_Xe_{\lambda,\theta}q\mu$ imply (4.4).

Proposition 4.1 is proved.

\vskip 2 mm
{\bf 5. Reconstruction of function $\psi_{\theta}\big|_{bX}$ from
Dirichlet-to-Neumann data on $bX$. Proof of Theorem 1.2A}

\vskip 2 mm
Let $X$ be domain containing $V_0$ and relativement compact in $V$ with
smooth (of classe $C^{(2)}$) boundry. Let $\sigma\in C^{(2)}(V)$, $\sigma>0$
on $V$, $\sigma=1$ on $V\b X$. Let $q={dd^c\sqrt{\sigma}\over \sqrt{\sigma}}$.
Let $u\in C(bX)$ and $\tilde u\in W^{1,\tilde p}(X)$, $\tilde p>2$,
be solution of the Dirichlet problem $d\sigma d^c\tilde u\big|_X=0$,
$\tilde u\big|_{bX}=u$, where $d^c=i(\bar\pa-\pa)$, $d=\bar\pa+\pa$.
 Let
$\tilde\psi=\sqrt{\sigma}\tilde u$ and $\psi=\sqrt{\sigma}u$. Then
$$dd^c\tilde\psi={dd^c\sqrt{\sigma}\over \sqrt{\sigma}}\tilde\psi=q\tilde\psi\
 \ {\rm on}\ \ X,\ \  \tilde\psi\big|_{bX}=\psi.\eqno(5.1)$$
Let $\psi_0$ be solution of Dirichlet problem
$$dd^c\psi_0\big|_X=0,\ \ \psi_0\big|_{bX}=\psi\big|_{bX}.$$
Let
$$\hat\Phi\psi=\bar\pa\tilde\psi\big|_{bX}\ \ {\rm and}\ \
\hat\Phi_0\psi=\bar\pa\tilde\psi_0\big|_{bX}.\eqno(5.2)$$
Operator $\psi\big|_{bX}\mapsto\bar\pa\tilde\psi\big|_{bX}$ is equivalent to
the Dirichlet-to-Neumann operator

\noindent
$u\big|_{bX}\mapsto\sigma d^c\tilde u\big|_{bX}$.

\vskip 2 mm
{\bf Proposition 5.1.}
{\sl
Let $\psi=e^{\lambda(z_1+\theta z_2)}\mu$ be the Faddeev type function
associated with potential $q={dd^c\sqrt{\sigma}\over \sqrt{\sigma}}$
(see Definition 1.4), generic divisor $\{a_1,\ldots,a_g\}$
with support in $V\b\bar X$ and generic $\theta\in\C$. Then\
$\forall\lambda\in\C\b E_{\theta}:\ |\lambda|\ge
const(V,\{a_j\}\theta,\sigma)$ the
restriction $\psi\big|_{bX}$ of $\psi$ on $bX$ can be found from
Dirichlet-to-Neumann operator $\psi\big|_{bX}\to\sigma d^c\psi\big|_{bX}$
through the uniquely solvable Fredholm integral equation
$$\mu_{\theta}(z,\lambda)\big|_{bX}+
\int\limits_{\xi\in bX}g_{\lambda,\theta}(z,\zeta)
m_{-\lambda}(\hat\Phi-\hat\Phi_0)m_{\lambda}\mu_{\theta}(\zeta,\lambda)=
1+i\sum_{j=1}^gC_{j,\theta}(\lambda)g_{\lambda,\theta}(z,a_j),\eqno(5.3)$$
$$i\sum_{j=1}^g(a_{j,1}+\theta a_{j,2})^{-k}C_{j,\theta}(\lambda)+
\int\limits_{z\in bX}(z_1+\theta z_2)^{-k}
(\pa+\lambda(dz_1+\theta dz_2))\mu=0,\ \ k=2,\ldots,g+1,$$
where $g_{\lambda,\theta}(z,\xi)$- kernel of operator
$R_{\lambda,\theta}\circ\hat R_{\theta}$,
$$m_{-\lambda}(\hat\Phi-\hat\Phi_0)m_{\lambda}\mu_{\theta}(\zeta,\lambda)=
\int\limits_{w\in bX}e^{-\lambda(\zeta_1+\theta\zeta_2)}
(\Phi(\zeta,w)-\Phi_0(\zeta,w))e^{\lambda(w_1+\theta w_2)}
\mu_{\theta}(w,\lambda),\eqno(5.4)$$
$\Phi(\zeta,w)$, $\Phi_0(\zeta,w)$ are kernels of operators
$\hat\Phi$ and $\hat\Phi_0$, $m_{\pm\lambda}$ denote the multiplication
operators by $e^{\pm\lambda(z_1+\theta z_2)}$,
values $\{a_{j,1}\}$ of the first coordinate
of points $\{a_j\}$ are supposed to be mutually different.
}

This proposition for the case $V=\C$ is equivalent to the second part of
Theorem 1 from [N2].

\vskip 2 mm
{\bf Lemma 5.1.}
{\sl
Let $\psi=e^{\lambda(z_1+\theta z_2)}\mu$ be Faddeev type function
of Proposition 5.1. Then\ $\forall\ z\in V\b X$ and
$\forall\lambda\in\C\b E_{\theta}:\ |\lambda|\ge
const(V,\{a_j\},\theta,\sigma)$
we have equalities
$$\eqalign{
&\mu_{\theta}(z,\lambda)=1-\int\limits_{\xi\in bX}g_{\lambda,\theta}(z,\xi)
\bar\pa\mu_{\theta}(\xi,\lambda)-\int\limits_{\xi\in bX}\mu_{\theta}(z,\xi)
e^{\lambda(\xi_1+\theta\xi_2)}\pa\bigl(e^{-\lambda(\xi_1+\theta\xi_2)}
g_{\lambda,\theta}(z,\xi)\bigr)+\cr
&i\sum_{j=1}^gC_{j,\theta}(\lambda)g_{j,\theta}(z,a_j)\cr}\eqno(5.5)$$
and
$$-\int\limits_{\xi\in bX}(z_1+\theta z_2)^{-k}
(\pa+\lambda(dz_1+\theta dz_2))\mu=
\sum_{j=1}^g(a_{j,1}+\theta a_{j,2})^{-k}iC_{j,\theta}(\lambda),\ k=2,\ldots
\eqno(5.6)$$
}

\vskip 2 mm
{\it Proof of Lemma 5.1.}

The equation
$$\bar\pa(\pa+\lambda(dz_1+\theta dz_2))\mu={i\over 2}q\mu+
i\sum_{j=1}^gC_{j,\theta}(\lambda)\delta(z,a_j),\eqno(5.7)$$
where $supp\,q\subseteq X$ implies that (1,0)-form
$f=(\pa+\lambda(dz_1+\theta dz_2))\mu$ is holomorphic on
$(V\b (X\cup_{j=1}^g\{a_j\})$ and
$Res_{a_j}(\pa+\lambda(dz_1+\theta dz_2))\mu={iC_j\over 2\pi i}$.
This and the property (4.12)
imply that\ $\forall\lambda\in\C\b E_{\theta}$ and\
$\forall\ k\ge 2$ form $(z_1+\theta z_2)^{-k}f$ is holomorphic in the
neighborhood of $(\tilde V\b V)$. By residue theorem applied to the form
$(z_1+\theta z_2)^{-k}f$ on $\tilde V\b X$, we obtain
$$\int\limits_{z\in bX}(z_1+\theta z_2)^{-k}f(z,\lambda)=
-2\pi i\sum_{j=1}^gRes_{a_j}(z_1+\theta z_2)^{-k}f(z,\lambda)=-
(a_{j,1}+\theta a_{j,2})^{-k}(iC_{j,\theta}(\lambda)),$$
$k=2,3,\ldots.$
Equality (5.6) is proved.

Let us prove now (5.5). Differential equation (5.7),  where
$\mu\big|_Y\in L^{\tilde p}(Y)$,

\noindent
$\mu\big|_{V\b\bar Y}\in L^{\infty}(V\b\bar Y)$, $\mu(z)\to 1$, $z\to\infty$,
$z\in V_1$, is equivalent by Lemma 3.1
to the system of equations
$$\eqalignno{
&\mu(z,\lambda)=1+R_{\lambda,\theta}\circ\hat R_{\theta}\bigl({i\over 2}q\mu+
i\sum_{j=1}^gC_j\delta(z,a_j)\bigr),\ z\in V,\ \ {\rm and}&(5.8)\cr
&\bar\pa(\pa+\lambda(dz_1+\theta dz_2))\mu=0,\ \
z\in V\b (X\cup_{j=1}^g\{a_j\}).&(5.9)\cr}$$
These equations imply relations (5.6). Besides, we have
equality
$$\int\limits_Xg_{\lambda,\theta}(z,\xi){i\over 2}q(\xi)\mu(\xi)=
\int\limits_Xg_{\lambda,\theta}(z,\xi)
\bar\pa(\pa+\lambda(dz_1+\theta dz_2))\mu.$$
Using Green-Riemann formula we obtain
$$\eqalign{
&\int\limits_Xe^{\lambda((z_1-\xi_1)+\theta(z_2-\xi_2))}
g_{\lambda,\theta}(z,\xi)\pa\bar\pa\psi=
\int\limits_X\psi\pa\bar\pa\bigl(e^{\lambda((z_1-\xi_1)+\theta(z_2-\xi_2))}
g_{\lambda,\theta}(z,\xi)\bigr)+\cr
&\int\limits_{bX}e^{\lambda((z_1-\xi_1)+\theta(z_2-\xi_2))}
g_{\lambda,\theta}(z,\xi)\bar\pa\psi+
\int\limits_{bX}\psi\pa\bigl(e^{\lambda((z_1-\xi_1)+\theta(z_2-\xi_2))}
g_{\lambda,\theta}(z,\xi)\bigr).\cr}$$
For $z\in V\b X$ we have
$$\pa\bar\pa\bigl(e^{\lambda((z_1-\xi_1)+\theta(z_2-\xi_2))}
g_{\lambda,\theta}(z,\xi)\bigr)=0.$$
Then
$$-\int\limits_{\xi\in X}g_{\lambda,\theta}(z,\xi)\bigl({i\over 2}q\mu\bigr)=
\int\limits_{\xi\in bX}g_{\lambda,\theta}\bar\pa\mu+
\int\limits_{\xi\in bX}e^{\lambda(\xi_1+\theta\xi_2)}\mu
\pa\bigl(e^{-\lambda(\xi_1+\theta\xi_2)}g_{\lambda,\theta}(z,\xi)\bigr).
\eqno(5.10)$$
From (5.8), (5.10) we deduce statement (5.5) of Lemma 5.1.

\vskip 2 mm
{\it Proof of Proposition  5.1.}

Let $\psi_0:\ \bar\pa\pa\psi_0\big|_X=0$ and $\psi_0\big|_{bX}=\psi$.
By Green-Riemann formula \ $\forall\ z\in V\b X$ we have
$$\int\limits_{\xi\in bX}\psi\pa
\bigl(e^{\lambda((z_1-\xi_1)+\theta(z_2-\xi_2))}
g_{\lambda,\theta}(z,\xi)\bigr)+
\int\limits_{\xi\in bX}
e^{\lambda((z_1-\xi_1)+\theta(z_2-\xi_2))}g_{\lambda,\theta}(z,\xi)
\bar\pa\psi_0=0.\eqno(5.11)$$
Formulas (5.11) and (5.5), (5.6) imply
$$\eqalign{
&\psi(z,\lambda)=e^{\lambda(z_1+\theta z_2)}-
\int\limits_{bX}
e^{\lambda((z_1-\xi_1)+\theta(z_2-\xi_2))}g_{\lambda,\theta}(z,\xi)
(\bar\pa\psi(\xi)-\bar\pa\psi_0(\xi))+\cr
&i\sum_{j=1}^ge^{\lambda(z_1+\theta z_2)}C_jg_{\lambda,\theta}(z,a_j).\cr}
\eqno(5.12)$$
Formula (5.12), (5.6) are equivalent to (5.3). Integral equation (5.3) is the
Fredholm equation in $C(bX)$, because operator $(\hat\Phi-\hat\Phi_0)$ is
compact operator in $C(bX)$. Existence\

\noindent
$\forall\lambda\in\C\b E_{\theta}$ of
unique Faddeev type function $\psi=e^{\lambda(z_1+\theta z_2)}\mu$,
associated with $q$ and divisor $\{a_1,\ldots,a_g\}$ imply existence
of solution of (5.3) with residue data

\noindent
$iC_j=Res_{a_j}(\pa+\lambda(dz_1+\theta dz_2))\mu,\ \ j=1,\ldots,g$.
Let us prove uniqueness of solution (5.3) in $C(bX)$ with residue data
$\{C_j\}$. Suppose $\mu\in C(bX)$ solves (5.3), (5.6). Consider this $\mu$ as
Dirichlet data for equation
$\bar\pa(\pa+\lambda(dz_1+\theta dz_2))\mu={i\over 2}q\mu$ on $X$, solution
of which well defines $\mu$ on $\bar X$.

Let us also define $\mu$ on $V\b\bar X$ by (5.5). Function $\mu(z,\lambda)$
defined in such a way on $V$ belongs to $C(V\b \cup_{j=1}^g\{a_j\})$.

Let us show that $\mu$ satisfy (5.7).
By Sohotsky-Plemelj jump formula\ $\forall\ z^*\in bX$ we have
$$\eqalign{
&{i\over 2}\mu(z^*)=\lim_{\scriptstyle z\to z^* \atop z\in X}\bigl(
\int\limits_{bX}g_{\lambda,\theta}\bar\pa\mu+
\mu e^{\lambda(\xi_1+\theta\xi_2)}\pa\bigl(e^{-\lambda(\xi_1+\theta\xi_2)}
g_{\lambda,\theta}\bigr)\bigr)-\cr
&-\lim_{\scriptstyle z\to z^* \atop z\in V\b X}\bigl(
\int\limits_{bX}g_{\lambda,\theta}\bar\pa\mu+
\mu e^{\lambda(\xi_1+\theta\xi_2)}\pa\bigl(e^{-\lambda(\xi_1+\theta\xi_2)}
g_{\lambda,\theta}\bigr)\bigr).\cr}\eqno(5.13)$$
From (5.5) and (5.13) we deduce equality
$$
\mu-{i\over 2}\mu=1-
\int\limits_{\xi\in bX}g_{\lambda,\theta}\bar\pa\mu-
\int\limits_{\xi\in bX}\mu
e^{\lambda(\xi_1+\theta\xi_2)}\pa\bigl(e^{-\lambda(\xi_1+\theta\xi_2)}
g_{\lambda,\theta}\bigr)+i\sum_{j=1}^gC_jg_{\lambda,\theta}(z,a_j),\ z\in X.
\eqno(5.14)$$
By Green-Riemann formula we have also
$$\eqalign{
&-\int\limits_{bX}g_{\lambda,\theta}\bar\pa\mu-
\mu e^{\lambda(\xi_1+\theta\xi_2)}\pa\bigl(e^{-\lambda(\xi_1+\theta\xi_2)}
g_{\lambda,\theta}\bigr)+i\sum_{j=1}^gC_{j,\theta}g_{\lambda,\theta}(z,a_j)=
\cr
&-\int\limits_X\mu(\bar\pa(\pa+\lambda(dz_1+\theta dz_2)g_{\lambda,\theta})+
\int\limits_Xg_{\lambda,\theta}\bar\pa(\pa+\lambda(dz_1+\theta dz_2))\mu+
i\sum_{j=1}^gC_{j,\theta}g_{\lambda,\theta}(z,a_j)=\cr}$$
$$\left\{\matrix{
&{\mu\over 2i}+
\int\limits_Xg_{\lambda,\theta}\bar\pa(\pa+\lambda(dz_1+\theta dz_2))\mu
+i\sum\limits_{j=1}^gC_{j,\theta}g_{\lambda,\theta}(z,a_j),\ z\in X,\hfill\cr
&\int\limits_Xg_{\lambda,\theta}\bar\pa(\pa+\lambda(dz_1+\theta dz_2))\mu
+i\sum\limits_{j=1}^gC_{j,\theta}g_{\lambda,\theta}(z,a_j),\ z\in V\b(X\cup_{j=1}^g
\{a_j\}).\hfill\cr}\right.\eqno(5.15)$$
Equalities (5.5), (5.6), (5.14) and (5.15) imply  (3.3) and
$$\eqalign{
&\mu(z)=1+
\int\limits_Vg_{\lambda,\theta}\bar\pa(\pa+\lambda(dz_1+\theta dz_2))\mu
+i\sum_{j=1}^gC_{j,\theta}g_{\lambda,\theta}(z,a_j)=\cr
&1+R_{\lambda,\theta}\circ\hat R_{\theta}\bigl({i\over 2}q\mu+
i\sum_{j=1}^gC_{j,\theta}\delta(z,a_j)\bigr),\ z\in V.\cr}$$
By Lemma 3.1 function $\mu_{\theta}(z,\lambda)$ is the Faddeev type function
associated with $q$ and divisor $\{a_1,\ldots,a_g\}$.
The uniqueness of solution of (5.3) in $C(bX)$ with residue data
$\{C_{j,\theta}\}$ follows now from uniqueness of the Faddeev type function.

\vskip 2 mm
{\bf 6. Reconstruction of conductivity function from Dirichlet-to-Neumann
data. Proof of Theorem 1.2B}

We will obtain here exact formulas for reconstruction of conductivity
function

\noindent
$\sigma\in C^{(3)}(V)$, $\sigma>0$, $\sigma\equiv 1$ on $V\b X$,
from Dirichlet-to-Neumann data
$$\psi_{\theta}\big|_{bX}\to  \bar\pa\psi_{\theta}\big|_{bX}$$
for Faddeev type functions
$\psi_{\theta}(z,\lambda)=e^{\lambda(z_1+\theta z_2)}\mu_{\theta}(z,\lambda)$,
$\theta\in \C\b \{\theta_1,\theta_d\}$,

\noindent
$\lambda\in\C\b E_{\theta}:\ |\lambda|\ge const(V,\{a_j\},\theta,\sigma)$,
$\{a_1,\ldots,a_g\}\subset Y\b X$.

For simplicity of presentation we consider in detail  the case of regular
algebraic curves in $\C^2\subset {\C}P^2$, only.

Let
$\tilde V=\{\tilde z=(\tilde z_0:\tilde z_1:\tilde z_2)\in {\C}P^2:\
\tilde P(\tilde z)=0\}$, where $\tilde P(\tilde z)$ homogeneous polynomial
of degre $N$. Let  ${\C}P^1_{\infty}=\{\tilde z:{\C}P^2:\tilde z_0=0\}$.
Put
$$\eqalign{
&\C^2=\{\tilde z\in {\C}P^2:\tilde z_0\ne 0\},\ \
z_1={\tilde z_1\over \tilde z_0},\ \ z_2={\tilde z_2\over \tilde z_0},\ \
P(z)=\tilde P(1,z_1,z_2),\cr
&V=\{z\in\C^2:\ P(z)=0\}=\tilde V\cap\C^2.\cr}\eqno(6.1)$$
Without restriction of generality we suppose that $\tilde V$ be (regular)
curve of degree $N\ge 2$ with property:
$$\eqalignno{
&\tilde V\cap {\C}P^1_{\infty}=\{\beta_1,\ldots,\beta_d\},\ \ {\rm where}\ \
\beta_1,\ldots,\beta_d\ \ {\rm be\ different\ points\ of}\ \ {\C}P^1_{\infty},
\cr
&\beta_l=(0:\beta_l^1:\beta_l^2),\ \ {\beta_l^2\over \beta_l^1}\in\C,\ \
l=1,\ldots,N,\cr
&{\pa P\over \pa z_2}(z)\ne 0,\ \ {\rm if}\ \ z\in V:\ |z_1|\ge r_0=const(V).
&(6.2)\cr}$$

For $\theta\in\C$ let $\{w_m\}$ be points of $V$, where
$(dz_1+\theta dz_2)\big|_V(w_m)=0$. Then for almost all $\theta\in\C$
the following relations are valid
$$\eqalign{
&\theta={\pa P\over \pa z_2}(w_m)\big/ {\pa P\over \pa z_1}(w_m),\ \
{\pa P\over \pa z_1}(w_m)\ne 0,\cr
&\big[{\pa^2 P\over \pa z_1^2}\bigl({\pa P\over \pa z_2}\bigr)^2-
2{\pa^2 P\over \pa z_1\pa z_2}\bigl({\pa P\over \pa z_2}\bigr)
\bigl({\pa P\over \pa z_1}\bigr)
+{\pa^2 P\over \pa z_2^2}\bigl({\pa P\over \pa z_1}\bigr)^2\big](w_m)\ne 0.
\cr}$$
Without restriction of generality it is sufficient to give proof under
condition that $\theta=0$, i.e. for points $w_m=(w_{m,1},w_{m,2})\in V$ such
that
$${\pa P\over \pa z_1}(w_m)\ne 0,\ \ {\pa P\over \pa z_2}(w_m)=0,\ \
{\pa^2 P\over \pa z_2^2}(w_m)\ne 0 \eqno(6.3)$$
and also such that \ $\forall\ m$ the line $\{z\in\C^2:\ z_1=w_{m,1}\}$ has
tangency with $X$ only in the single point $w_m$, $m=1,\ldots,M$.
By  Hurwitz-Riemann formula $M=N(N-1)$. In the neighborhood of point
$w_m\in V$ curve $V$ can be represented in the form
$$\eqalign{
&V=\{(z_1,z_2)\in\C^2:\ z_1=w_{m,1}+\cr
&\bigl({\pa P\over \pa z_1}(w_m)\bigr)^{-1}\bigl[-{1\over 2}
{\pa^2 P\over \pa z_2^2}(w_m)(z_2-w_{m,2})^2+O((z_2-w_{m,2})^3)\bigr].\cr}
\eqno(6.4)$$
The reconstruction formula for ${dd^c\sqrt{\sigma}\over \sqrt{\sigma}}(w_m)$,
$m=1,\ldots,M$, will be obtained here by the stationary phase method, using
formula (4.17).

Let $\mu$ be Faddeev type function (3.1) with properties (3.1a)-(3.1c) and
with $\theta=0$.

Below in this section we will write $\hat R_0$, $R_{\lambda,0}$,
$e_{\lambda,0}$, $\mu_0$, $\psi_0$, $\Delta_0$, $E_0$, $C_{j,0}$ as
$\hat R$, $R_{\lambda}$, $e_{\lambda}$, $\mu$, $\psi$, $\Delta$, $E$, $C_j$.

Let
$$f_0=F_0dz_1={i\over 2}\hat R(q\mu),\ \
f_1=F_1dz_1=i\sum_{j=1}^gC_j(\lambda)\hat R(\delta(\cdot,a_j)),$$
where $\mu=\mu(z,\lambda)$, $z\in V$,
$\lambda\in\C\b E:\ |\lambda|\ge const(V,\{a_j\},\sigma)$.

\vskip 2 mm
{\bf Lemma 6.1.}
{\sl
For $u_0=R_{\lambda}f_0$ the following estimate holds:
$$\|u_0(\cdot,\lambda)-
{F_0(\cdot,\lambda)\over \lambda}\|_{L^{9/4}(X)}\le
{const(V,\tilde p)\over |\lambda|^{7/5}}
\|f_0(\cdot,\lambda)\|_{\tilde W_{1,0}^{2,\tilde p}(V)}.$$
}

\vskip 2 mm
{\it Proof of Lemma 6.1.}
By Lemma 2.1 and Proposition 2 from [He] we have
$f_0\in \tilde W_{1,0}^{2,\tilde p}(V)$, $F_0\in \tilde W^{1,p}(V)$.
Using equality $\pa_ze_{\lambda}(z)=\lambda e_{\lambda}(z)dz_1$ and
integration by parts formula $u_0=R_{\lambda}f_0=e_{-\lambda}(z)
\overline{R(\overline{e_{\lambda}f_0})}$ can be transformed into the
following
$$\eqalign{
&u_0(z)=e_{-\lambda}(z)\overline{R_1(\overline{e_{\lambda}f_0})}+
e_{-\lambda}(z)\overline{R_0(\overline{e_{\lambda}f_0})}=\cr
&-{e_{-\lambda}(z)\over 2\pi i}{1\over \lambda}\int\limits_V
{e_{\lambda}(\xi)\pa F_0\wedge\overline{d\xi_1}
\det\bigl[{\pa\bar P\over \pa\bar\xi}(\xi),\xi-z\bigr]\over
{\pa\bar P\over \pa\bar\xi_2}(\xi)\cdot |\xi-z|^2}-\cr
&-{e_{-\lambda}(z)\over 2\pi i}{1\over \lambda}\int\limits_V
e_{\lambda}(\xi)F_0\pa\biggl(
{\det\bigl[{\pa\bar P\over \pa\bar\xi}(\xi),\xi-z\bigr]\wedge d\bar\xi_1\over
{\pa\bar P\over \pa\bar\xi_2}(\xi)\cdot |\xi-z|^2}\biggr)+
e_{-\lambda}(z)\overline{R_0(\overline{e_{\lambda}f_0})},\cr}\eqno(6.5)$$
where $R_1$, $R_0$ operators defined in section 1 (see remark 1.1).
From (6.5), using

\noindent
Corollary 1.2 from [He], we deduce
$$\lambda u_0-F_0=-e_{-\lambda}(z)
\overline{R_1(\overline{e_{\lambda}(\xi)\pa F_0})}-
e_{-\lambda}(z)\overline{R_0(\overline{e_{\lambda}(\xi)\pa F_0})}\buildrel
\rm def \over =J_1(z)+J_0(z).\eqno(6.6)$$
We will estimate further only term $J_1(z)$. Estimate for $J_0(z)$ is
similar.

For $J_1(z)$ we have $J_1(z)=J_1^+(z)+J_1^-(z)$, where
$$J_1^{\pm}(z)=
{e_{-\lambda}(z)\over 2\pi i}\int\limits_V
{e_{\lambda}(\xi)\chi_{\rho}^{\pm}(\xi)\pa F_0(\xi)\wedge\overline{d\xi_2}
\det\bigl[{\pa\bar P\over \pa\bar\xi}(\xi),\xi-z\bigr]\over
{\pa\bar P\over \pa\bar\xi_1}(\xi)\cdot |\xi-z|^2},$$
$$\eqalign{
&\chi_{\rho}^{\pm}\ \ {\rm be\ smooth\ functions\ such\ that}\ \
\chi_{\rho}^++\chi_{\rho}^-\equiv 1,\cr
&\chi_{\rho}^+=1,\ \ {\rm if}\ \ \big|{d\xi_1\over d\xi_2}\big|\le\rho,\ \
supp\,\chi_{\rho}^+\subset \{\xi:\ \big|{d\xi_1\over d\xi_2}\big|\le 2\rho\}
\cr
&{\rm and}\ \ |d\chi_{\rho}^{\pm}|=O({1\over \rho}).\cr}\eqno(6.7)$$
Let $B_0=\{z\in V:\ d\xi_1\big|_V(z)=0\}$. Property
$\bar\pa F_0=dz_1\rfloor {i\over 2}q\mu$ implies estimate

\noindent
$\bar\pa F_0=O\bigl({1\over dist(z,B_0)}\bigr)dz_2$. From this, formula for
$J_1^+(z)$ and Lemma 3.1 of [He] we obtain estimate for
$$J_1^+:\ \|J_1^+\|_{L^{9/4}(X)}=O(\rho^{2/3})
\|f_0\|_{\tilde W^{2,\tilde p}_{1,0}(V)}.\eqno(6.8)$$

In order to estimate $J_1^-(z)$ we  integrate by parts in the formula for
$J_1^-$, using

\noindent
$\pa_ze_{\lambda}(z)=\lambda e_{\lambda}(z)dz_1$. Then
inequalities
$$|\bar\pa F_0(z)|=O\bigl({1\over dist(z,B_0)}\bigr),\ \
|\bar\pa\pa F_0(z)|=O\bigl({1\over dist(z,B_0)}\bigr)^2,\ \ z\in X\b B_0$$
and inequality
$$\bigg|\int\limits_{\rho\le |\xi_2|\le 1}
{d\xi_2\wedge d\bar\xi_2\over |\xi_2|^2(\bar\xi_2-\bar z_2)}\bigg|+
\bigg|\int\limits_{\rho\le |\xi_2|\le 1}
{d\xi_2\wedge d\bar\xi_2\over |\xi_2|(\bar\xi_2-\bar z_2)^2}\bigg|=
O({1\over \rho})$$
imply estimate
$$\|J_1^-\|_{L^{\infty}(X)}=O({1\over |\lambda|\rho})
\|f_0\|_{\tilde W_{1,0}^{2,\tilde p}(V)}.\eqno(6.9)$$
From (6.6), (6.8), (6.9) with $\rho=|\lambda|^{-3/5}$
we obtain statement of Lemma 6.1.

\vskip 2 mm
{\bf Lemma 6.2.}
{\sl
Let $q\in C_{1,1}^{(1)}(V)$, $supp\,q\subseteq X$,
$f_0={i\over 2}\hat R(q\mu)$, $u_0=R_{\lambda}f_0$. Then the following
asymptotic estimate is valid
$$\eqalign{
&\big|\int\limits_Xe_{\lambda}(z)q(z)u_0(z)\big|=o({1\over |\lambda|}),\ \
{\it for}\cr
&\lambda\in\C:\ |\lambda|\ge const(V,\{a_j\},\sigma),\ \
|\Delta(\lambda)(1+|\lambda|)^g|\ge\delta>0,\ \ {\it for\ some\ sufficient\
small}\ \ \delta.\cr}$$
}

\vskip 2 mm
{\bf Proof of Lemma 6.2.}

From Lemma 6.1, using estimate of $\mu$ from (3.1b), we obtain
asymptotic relation in the space $L^{\tilde p}(V)$, $2<\tilde p<9/4$:
$$\eqalign{
&u_0(z,\lambda)={F_0(z,\lambda)\over \lambda}+O({1\over |\lambda|^{7/5}})=\cr
&{dz_1\rfloor {i\over 2}\hat R(q)\over \lambda}+O({1\over |\lambda|^{7/5}})\
\ {\rm if}\ \ |\Delta_{\theta}(\lambda)(1+|\lambda|)^g|\ge\delta>0.\cr}$$
Putting this relation into $\int\limits_Xe_{\lambda}(z)q(z)u_0(z)$, we obtain
$$\int\limits_Xe_{\lambda}(z)q(z)u_0(z)={i\over 2\lambda}
\int\limits_Xe_{\lambda}(z)q(z)(dz_1\rfloor\hat R(q))+
O({1\over |\lambda|^{7/5}}).$$
By Riemann-Lebesgue type theorem
$$\int\limits_Xe_{\lambda}(z)q(z)(dz_1\rfloor\hat R(q))=o(1)\ \ {\rm if}\ \
\lambda\to\infty,\ \ |\Delta(\lambda)|(1+|\lambda|)^g\ge\delta>0.$$
This implies the statement of Lemma 6.2.

\vskip 2 mm
{\bf Lemma 6.3.}
{\sl
Let $q\in C_{1,1}^{(1)}(V)$, $supp\,q\subset X$.
Let $w_1,\ldots,w_M$ be the points, where $dz_1\big|_V(w_m)=0$. Then the
following consequence of stationary phase method is valid:
$$\int\limits_Ve^{i\tau(z_1+\bar z_1)}q(z)=\sum_m(1+o(1))\sum_{m=1}^M
-{\pi\over r}{\big|{\pa P\over \pa z_1}(w_m)\big|Q_2(w_m)\over
\big|{\pa^2 P\over \pa z_2^2}(w_m)\big|}e^{i\tau(w_{m,1}+\bar w_{m,1})},
\eqno(6.10)$$
where $Q_2(w_m)={q\over 2idz_2\wedge d\bar z_2}(w_m)$.
}

\vskip 2 mm
{\it Proof of Lemma 6.3} (see [Fe], Th.2.1).

\vskip 2 mm
{\bf Lemma 6.4.}
{\sl
Let $q={dd^c\sqrt{\sigma}\over \sqrt{\sigma}}\in C_{1,1}^{(1)}(X)$,
$supp\,q\subset X$,
$f_1=i\sum\limits_{j=1}^gC_j(\lambda)\hat R(\delta(\cdot,a_j))$,
$u_1=R_{\lambda}f_1$. Then the following asymptotic estimate is valid
$$\eqalign{
&\big|\int\limits_Xe_{\lambda}(z)q(z)u_1(z)\big|=
O({1\over |\lambda|^{3/2-\ep}}),
\ \ {\it for}\ \
\lambda\in\C:\ |\lambda|\ge const(V,\{a_j\},\sigma,\ep),\cr
&|\Delta(\lambda)(1+|\lambda|)^g|\ge\delta>0,\ \ \delta\ \
{\it for\ some\ sufficiently\ small}\ \ \delta.\cr}$$
}

\vskip 2 mm
{\bf Proof of Lemma 6.4.}
Using that $\{a_1,\ldots,a_g\}$ be generic divisor, from  estimate (3.7)
(Lemma 3.3) we obtain inequality
$$\sup\limits_{j,\lambda}|C_j(\lambda)|\le const(V,\{a_j\})\sup\limits_k
\bigg|\int\limits_Xe_{\lambda}(z)
\big(i{dd^c\sqrt{\sigma}\over \sqrt{\sigma}}+2\bar\pa\ln\sqrt{\sigma}\wedge
\pa\ln\sqrt{\sigma}\bigr)\mu{\bar\omega_k\over d\bar z_1}(z)\bigg|.$$
Let $\ep>0$ be small enough and $B_{\ep}=
\{z\in X:\ \big|{dz_1\over dz_2}\big|<\ep\}$. Then
$$\big|{\bar\omega_k(z)\over d\bar z_1}\big|_X=O\bigl(\sum_{m=1}^M
{1\over |z_2-w_{m,k}|}\bigr),\ \ z\in X.$$

Let $\chi_{\rho}^{\pm}\in C^{(1)}(X)$ be functions with properties (6.7).
Using that $\sigma\in C^{(3)}(X)$, $\mu\in\tilde W^{1,\tilde p}(X)$,
$\pa_ze_{\lambda}(z)=\lambda e_{\lambda}(z)dz_1$ by integration by parts we
obtain
$$\eqalign{
&\bigg|\int\limits_X\chi_{\rho}^-(z)e_{\lambda}(z)
\big(i{dd^c\sqrt{\sigma}\over \sqrt{\sigma}}+2\bar\pa\ln\sqrt{\sigma}\wedge
\pa\ln\sqrt{\sigma}\bigr)\mu(z,\lambda){\bar\omega_k\over d\bar z_1}(z)\bigg|
\le\cr
&{const(V,\sigma)\over \rho \lambda}.\cr}$$
We have also directly
$$\eqalign{
&\bigg|\int\limits_X\chi_{\rho}^+(z)e_{\lambda}(z)
\big(i{dd^c\sqrt{\sigma}\over \sqrt{\sigma}}+2\bar\pa\ln\sqrt{\sigma}\wedge
\pa\ln\sqrt{\sigma}\bigr)\mu(z,\lambda){\bar\omega_k\over d\bar z_1}(z)\bigg|
\le\cr
&const(V,\sigma)\rho.\cr}$$
These estimates with $\rho={1\over \sqrt{|\lambda|}}$ and estimates for
Faddeev type Green  function

\noindent
$|R_{\lambda}\circ\hat R(\delta(\cdot,a_j)|=O({1\over |\lambda|^{1-\ep}})$
from Theorem 4 of [He] imply statement of Lemma 6.4.

\vskip 2 mm
{\bf Proposition 6.1.}
{\sl
Under conditions (6.1)-(6.4), for $\lambda=i\tau:\ \tau\in\R$,
$|\tau^g\Delta(i\tau)|\ge\delta>0$, $\delta$- small enough, the
following formula is valid
$$\eqalign{
&\int\limits_{z\in bX}e_{i\tau}(z)\bar\pa_z\mu(z,i\tau)=
\int\limits_{z\in X}e_{i\tau}(z){q\mu\over 2i}=\cr
&{{1+o(1)}\over \tau}\sum_{m=1}^M{\pi i\over 2}
{dd^c\sqrt{\sigma}\over \sqrt{\sigma}dd^c|z|^2}\big|_V(w_m)
e^{i\tau(w_{m,1}+\bar w_{m,1})}
\big|{\pa^2 P\over \pa z_2^2}(w_m)\big|^{-1}
{\pa P\over \pa z_1}(w_m).\cr}\eqno(6.11)$$
}

\vskip 2 mm
{\it Proof of Proposition 6.1 and Theorem 1.2B}.
From Lemma 3.1 we have equality
$$\mu=1+R_{\lambda}\circ\hat R\bigl({i\over 2}q\mu\bigr)+
R_{\lambda}\circ\hat R\bigl(i\sum_{j=1}^gC_j\delta(z,a_j)\bigr)=1+u_0+u_1.
\eqno(6.12)$$

Let $\delta>0$ be small enough.
Estimates of Lemmas 6.2, 6.4 and (6.12) give asymptotic equality
$$\mu=1+o\bigl({1\over \lambda}\bigr) \eqno(6.13)$$
under conditions $\lambda\in\C:\ |\lambda|\ge const(V,\{a_j\},\sigma)$,
$|\Delta(\lambda)(1+|\lambda|)^g|\ge\delta>0$.

By Proposition 1.1,\ $\forall\ep>0$ we have inequality
$$\underline{\lim}_{\lambda\to\infty}
|\lambda^g\Delta(\lambda)|_{\ep}=\delta(\ep)>0,\ \ {\rm where}\ \
 |\lambda^g\Delta(\lambda)|_{\ep}=
\sup\limits_{\{\lambda^{\prime}:\ |\lambda^{\prime}-\lambda|\le\ep\}}
|\lambda^{\prime}\Delta(\lambda^{\prime}|.$$
So for any $\ep>0$ and any positive $\delta<\delta(\ep)$ there exists $r$
such that the set

\noindent
$\{\lambda\in\C:\ |\Delta(\lambda)(1+|\lambda|)^g|\ge\delta>0\}$
intersects any disque $\{\lambda^{\prime}:\ |\lambda-\lambda^{\prime}|<\ep\}$, with
$|\lambda|\ge r.$
This property, Lemma 6.3 and property (6.13) imply Proposition 6.1.

Theorem 1.2B follows from Proposition 1. Indeed, stationary phase method
permits differentiation of (6.11) with respect to $\tau$, keeping
(in our case)
terms of order ${1\over \tau}$. Differentiation of the right-hand side of
(6.11) gives for $\theta=0$ the right-hand side of (1.12).

Theorem 1.2B is proved.

\vskip 2 mm
{\bf Remark 6.1.}
To obtain version of Proposition 6.1 with arbitrary generic $\theta$ from
Proposition 6.1 with $\theta=0$ it is sufficient to change coordinate
system:\ $\tilde z_1=z_1+\theta z_2$, $\tilde z_2=z_2$.

\vskip 2 mm
{\bf Remark 6.2.}
Proposition 6.1 can be reformulated also as formula for reconstruction of
conductivity function from scattering data $b_{\theta}(i\tau)$ and
$C_{j,\theta}(i\tau)$. Indeed,by formula (4.16), we have
$$\int\limits_{bX}e_{i\tau,\theta}(z)\bar\pa\mu(z,i\tau)=-2\pi
\bigl[\tau b_{\theta}(i\tau)d+\sum_{j=1}^gC_{j,\tau}(i\tau)
e_{i\tau,\theta}(a_j),$$
where $d$ is defined in section 1.

\vskip 2 mm
{\bf 7. Proof of Proposition 1.1}

For simplicity of presentation we give proof only for the case when $V$ is
algebraic curve in $\C^2$.
Proposition 1.1 will be obtained here as a corollary of the following
statement.

\vskip 2 mm
{\bf Proposition 7.1.}
{\sl
Let $\theta\in\C\b\{\theta_1,\ldots,\theta_N\}$,
$\delta=\delta(\theta)=\inf\limits_l|\theta-\theta_l|$,
$V_0=\{z\in V:\ |z_1|\le r_0(\delta)\}$,
$r_0(\delta)={const(V)\over \sqrt{\delta}}$.
Let $\{b_m\}$ be the points of $V$, where $(dz_1+\theta dz_2)\big|_V(b_m)=0$,
$m=1,\ldots,M$, and $\{a_1,\ldots,a_g\}$ be the points of generic divisor in
$V\b \bar V_0$.  Then \ $\forall\ j,k=1,\ldots,g$ and for $\lambda=i\tau$,
where $\tau\in\R$, large enough, such that
$|\Delta_{\theta}(i\tau)|\ge\delta>0$,
the following asymptotic equality is valid
$$\eqalign{
&\int\limits_V
\hat R_{\theta}(\delta(\xi,a_j))\wedge\bar\omega_k(\xi)
e_{\lambda,\theta}(\xi)=
-{1\over \bar\lambda}e_{\lambda,\theta}(a_j)
{\bar\omega_k\over {d\bar\xi_1+\bar\theta d\bar\xi_2}}(a_j)-\cr
&-{\pi\over |\lambda|}\sum_{m=1}^M
\exp{[\lambda(b_{m,1}+\theta b_{m,2})-\bar\lambda
(\bar b_{m,1}+\bar\theta\bar b_{m,2})]}K_{j,k}(b_m,a_j)+
O\bigl({1\over |\lambda|^2}\bigr),\cr}$$
where
$$\eqalign{
&K_{j,k}(b_m,a_j)=\cr
&{|{\pa P\over \pa z_1}(b_m)|^3\hat R_{\theta}(\delta(b_m,a_j))\wedge
\bar\omega_k(b_m)(1+|\theta|^2)\over
|{\pa^2P\over \pa z_1^2}\bigl({\pa P\over \pa z_2}\bigr)^2-2
{\pa^2P\over \pa z_1\pa z_2}{\pa P\over \pa z_2}{\pa P\over \pa z_1}+
{\pa^2P\over \pa z_2^2}\bigl({\pa P\over \pa z_1}\bigr)^2|
d d^c|z|^2\big|_V(b_m)}.\cr}\eqno(7.1)$$
}

\vskip 2 mm
{\bf Lemma 7.1.}
{\sl
Let $V\b V_0=\cup_{l=1}^gV_l$ be a curve with properties i)-iv) of section 1.
Then\ $\forall\theta\ne\theta_1,\ldots,\theta_d$  any point $w$, where
$(dz_1+\theta dz_2)\big|_V(w)=0$, belongs to
$V_0=\{z\in V:\ |z_1|\le r_0(\delta)\}$, where
$r_0(\delta)=const(V)/\sqrt{\delta}$, $\delta=\min\limits_l|\theta-\theta_l|$.
}

\vskip 2 mm
{\it Proof of Lemma 7.1.}
For any point $w\in V\b V_0$, where $(dz_1+\theta dz_2)\big|_V(w)=0$,
definition $\theta_l=-{1\over \gamma_l}$,

\noindent
$l=1,\ldots,d$, and
property iii) of Section 1 imply for some $l=l(w)$ equality
$$\eqalign{
&0=(dz_1+\theta dz_2)\big|_V(w)=
\bigl[1+\theta\bigl(\gamma_l+{\gamma_l^0\over w_1^2}+
O\bigl({1\over w_1^3}\bigr)\bigr)\bigr]dz_1=\cr
&\bigl[1+\theta\gamma_l+O\bigl({\theta\over w_1^2}\bigr)\bigr]dz_1=\gamma_l
\bigl[(\theta-\theta_l)+O\bigl({\theta\over w_1^2}\bigr)\bigr]dz_1.\cr}$$
This gives equality
$\theta\bigl(1+O\bigl({1\over w_1^2}\bigr)\bigr)=\theta_l$.
This equality together with inequality $|\theta-\theta_l|\ge\delta$ implies
inequality $|w_1|\le {const(V)\over \sqrt{\delta}}=r_0(\delta)$.

Lemma 7.1 is proved.

Let further
$$\eqalign{
&A_{\ep,j}=\{z\in V:\ |z-a_j|\le\ep\},\ A_{\ep}=\cup_{j=1}^gA_{\ep,j},\cr
&B_{\ep,m}=\{z\in V:\ |z-b_m|\le\ep\},\ B_{\ep}=\cup_{m=1}^MB_{\ep,m}.\cr}$$

\vskip 2 mm
{\bf Lemma 7.2}
{\sl
Let $r_0(\delta)$, $\delta=\delta(\theta)$ be as in Lemma 7.1. Let
$\chi^{A_{\ep}}$, $\chi^{B_{\ep}}$ be smooth functions with properties
$$\eqalign{
&\chi^{A_{\ep}}\big|_{A_{\ep}}=1,\ \ \chi^{A_{\ep}}\big|_{V\b A_{2\ep}}=0,\ \
|d\chi^{A_{\ep}}|=O({1\over \ep}),\cr
&\chi^{B_{\ep}}\big|_{B_{\ep}}=1,\ \ \chi^{B_{\ep}}\big|_{V\b B_{2\ep}}=0,\ \
|d\chi^{B_{\ep}}|=O({1\over \ep}).\cr}$$
Then for any $\ep>0$ small enough we have
$B_{2\ep}\cap A_{2\ep}=\{\emptyset\}$ and\  $\forall\ j,k=1,\ldots,g$
$$\Delta_{\theta,\ep}^{j,k}(\lambda)\buildrel \rm def \over =
\int\limits_{\xi\in V}(1-\chi^{A_{\ep}}-\chi^{B_{\ep}})
\hat R(\delta(\xi,a_j))\wedge\bar\omega_k(\xi)e_{\lambda,\theta}(\xi)=
O\bigl({1\over \lambda^2}\bigr).$$
}

\vskip 2 mm
{\it Proof of Lemma 7.2.}
By Lemma 7.1, any point $b_m$, where $(dz_1+\theta dz_2)\big|_V(b_m)=0$
belongs to
$\{z\in V:\ |z_1|\le r_0\}$. Under the  conditions of Lemma 7.2, any
$a_j$ from  $\{a_1,\ldots,a_g\}$ belongs to
$\{z\in V:\ |z_1|>r_0(\delta)\}$, $\delta=\delta(\theta)$.

Then $B_{2\ep}\cap A_{2\ep}=\{\emptyset\}$ if $\ep$ is  small enough.
From definition of $\Delta^{j,k}_{\theta,\ep}$ and equality
$\bar\pa\hat R_{\theta}(\delta(\ep,a_j))\big|_{V\b \{a_j\}}=0$ we obtain
$$\eqalign{
&\Delta^{j,k}_{\theta,\ep}(\lambda)={1\over \bar\lambda}
\int\limits_V(1-\chi^{A_{\ep}}-\chi^{B_{\ep}})
\hat R_{\theta}(\delta(\xi,a_j))
\wedge{\bar\omega_k\over {d\bar\xi_1+\bar\theta d\bar\xi_2}}
\bar\pa e_{\lambda,\theta}(\xi)=\cr
&-{1\over \bar\lambda}
\int\limits_V(1-\chi^{A_{\ep}}-\chi^{B_{\ep}})
\hat R_{\theta}(\delta(\xi,a_j))
\wedge\bar\pa\bigl({\bar\omega_k\over {d\bar\xi_1+\bar\theta d\bar\xi_2}}
\bigr)e_{\lambda,\theta}(\xi)-\cr
&{1\over \bar\lambda}\int\limits_V
\bar\pa(\chi^{A_{\ep}}+\chi^{B_{\ep}})
\hat R_{\theta}(\delta(\xi,a_j))
\wedge{\bar\omega_k(\xi)\over {d\bar\xi_1+\bar\theta d\bar\xi_2}}
e_{\lambda,\theta}(\xi)+\cr
&{1\over \bar\lambda}
\lim\limits_{r\to\infty}
\int\limits_{\{\xi\in V:\ |\xi_1|=r\}}
\hat R_{\theta}(\delta(\xi,a_j))
\wedge{\bar\omega_k\over {d\bar\xi_1+\bar\theta d\bar\xi_2}}
e_{\lambda,\theta}(\xi).\cr}\eqno(7.2)$$
From asymptotic estimates
$|\hat R_{\theta}(\delta(\xi,a_j))|=O(|d\xi_1|)$ and
$|\bar\omega_k|=O\bigl({d\bar\xi_1\over \bar\xi_1^2}\bigr)$, $\xi_1\to\infty$,
and property $\inf\limits_l|\theta-\theta_l|>0$ we obtain vanishing of the
last term of the right-hand side of (7.2).

Property $(d\xi_1+\theta d\xi_2)\big|_{V\b B_{\ep}}\ne 0$ permits to
integrate other terms of the right-hand side of (7.2) by parts once more and
to obtain statement of Lemma 7.2.

\vskip 2 mm
{\bf Lemma 7.3}
{\sl
For any $k,j\in\{1,\ldots,g\}$, $\theta\notin\{\theta_1,\ldots,\theta_d\}$
and any $\ep>0$ we have the asymptotic equality
$$\int\limits_V\chi^{A_{\ep,j}}
\hat R_{\theta}(\delta(\xi,a_j))\wedge\bar\omega_k(\xi)e_{\lambda,\theta}(\xi)
=-{1\over \bar\lambda}e_{\lambda,\theta}(a_j)
{\bar\omega_k\over {d\bar\xi_1+\bar\theta d\bar\xi_2}}(a_j)+
\bigl({1\over \lambda^2}\bigr).$$
}

\vskip 2 mm
{\it Proof of Lemma 7.3.}
Integration by parts of the left-hand side, equality

\noindent
$\bar\pa\hat R(\delta(\xi,a_j))=\delta(\xi,a_j)$ and inequality
$(d\xi_1+\theta d\xi_2)\big|_{A_{\ep,j}}\ne 0$ imply statement of
Lemma 7.3.

\vskip 2 mm
{\bf Lemma 7.4}
{\sl
Under the conditions of Lemmas 7.1, 7.2,\ $\forall\delta>0$,
$\theta:\ \inf\limits_l|\theta-\theta_l|>\delta$, $\forall\ j,k=1,\ldots,g$,
$$\eqalign{
&\int\limits_V\chi^{B_{\ep}}
\hat R_{\theta}(\delta(\xi,a_j))\wedge\bar\omega_k(\xi)e_{i\tau,\theta}(\xi)
=\cr
&-{\pi\over |\lambda|}\sum_{m=1}^M
\exp[\lambda(b_{m,1}+\theta b_{m,2})-\bar\lambda
(\bar b_{m,1}+\bar\theta\bar b_{m,2})]K_{j,k}(b_m,a_j)+
O\bigl({1\over |\lambda|^2}\bigr),\cr}$$
where $\theta=\theta(b_m)$, $m=1,\ldots,M$, and $K_{j,k}(b_m,a_j)$ are defined
 by (7.1).
}

\vskip 2 mm
{\it Proof of Lemma 7.4.}
This statement is consequence of the classical result of the stationary
phase method [Fe], applied to the left-hand side, taking into account the
following equality for $e_{\lambda,\theta}(z)$
in the neighborhood of the stationary points $b_m\in V$, $m=1,\ldots,M$,
$$e_{\lambda,\theta}(z)=
\exp[\lambda(b_{m,1}+\theta b_{m,2})-\bar\lambda
(\bar b_{m,1}+\bar\theta\bar b_{m,2})]\times
\exp[\lambda A(z_2-b_{m,2})^2-\bar\lambda\bar A(\bar z_2-\bar b_{m,2})^2],$$
where
$$A=-{
\bigl({\pa^2P\over \pa z_1^2}\theta^2-2
{\pa^2P\over \pa z_1\pa z_2}\theta+
{\pa^2P\over \pa z_2^2}\bigr)(b_m)(z_2-b_{m,2})^2(1+O(z_2-b_{m,2}))\over
2\bigl({\pa P\over \pa z_1}\bigr)(b_m)}.$$
We use here $z_2,\bar z_2$ as coordinates of integration.

Lemma 7.4 is proved.

\vskip 2 mm
{\it Proof of Proposition 7.1.}
Proposition 7.1 follows from Lemmas 7.2-7.4.

In the proof of Proposition 1.1 we will apply also the following statement
about exponential polynomials discovered by L.Ehrenpreis [E] and reinforced
by C.Berenstein and M.Dostal [BD].

\vskip 2 mm
{\bf Proposition 7.2.}  ([E], [BD])
{\sl
Let $Q(\xi)$ be an exponential polynomial
$$Q(\xi)=\sum_{k=1}^Nq_k(\xi)e^{<\alpha_k,\xi>},$$
where $\{q_k\}$ are polynomials of $\xi=(\xi_1,\ldots,\xi_n)\in\C^n$,
$\alpha_k=\{\alpha_{k,1},\ldots,\alpha_{k,n}\}\in\C^n$,

\noindent
$k=1,\ldots,N$.

Let $h(\xi)=\max\limits_kRe<\alpha_k,\xi>$. Then\ $\forall\ep>0\ \exists$
 constant $C=C(\ep,Q)>0$ such that
$$|Q(\xi)|_{\ep}\buildrel \rm def \over =
\sup\limits_{\{\xi^{\prime}\in\C:\ |\xi^{\prime}-\xi|<\ep\}}
|Q(\xi^{\prime})|\ge {1\over C}e^{h(\xi)}.$$
}

\vskip 2 mm
The final part of the proof of Proposition 1.1 consists of the following.

\noindent
Proposition 7.1  and definition of $\Delta_{\theta}(\lambda)$ imply
asymptotic equality
$$\eqalign{
&|\lambda|^g\Delta_{\theta}(\lambda)=\det\bigl(-{\lambda\over \bar\lambda}
e_{\lambda,\theta}{\bar\omega_k\over {d\bar\xi_1+\bar\theta d\bar\xi_2}}(a_j)
-\cr
&\pi\sum_{m=1}^M\exp[\lambda(b_{m,1}+\theta b_{m,2})-
\bar\lambda(\bar b_{m,1}+\bar\theta\bar b_{m,2})]K_{j,k}(b_m,a_j)\bigr)
+O\bigl({1\over |\lambda|}\bigr),\cr}\eqno(7.3)$$
where $j,k=1,\ldots,g$.

The determinant of the right-hand side of (7.3) is
an exponential polynomial  $Q(\lambda,\bar\lambda)$ of the form
$$Q(\lambda,\bar\lambda)=\sum_{k=1}^Nq_k(\lambda,\bar\lambda)
e^{\lambda\alpha_k-\bar\lambda\bar\alpha_k},\eqno(7.4)$$
where $\lambda\in\C$, $\alpha_k\in\C$, $k=1,\ldots,N$.
Coefficient $q_k(\lambda,\bar\lambda)$ of exponential polynomial
$Q(\lambda,\bar\lambda)$ and complex frequences $\{\alpha_k\}$ depend on
$V$, $\{a_j\}$, $\theta$, $\{b_m\}$.
Applying Proposition 7.2 to the exponential polynomial (7.4) we obtain
uniformly for $\lambda\in\C$ estimate
$$|Q(\lambda,\bar\lambda)\big|_{\ep}\ge {1\over C(\ep,Q)}
e^{\max\limits_kRe\,(\lambda\alpha_k-\bar\lambda\bar\alpha_k)}=
{1\over C(\ep,Q)}.\eqno(7.5)$$
The both inequalities of Proposition 1.1  follow
from (7.3)-(7.5).

\vskip 2 mm
{\bf References}

\item{[ BC1]} Beals R., Coifmann R., Multidimensional inverse scattering and
nonlinear partial differential equations, Proc. Symp.Pure Math. {\bf 43}
(1985), A.M.S. Providence, Rhode Island, 45-70
\item{[ BC2]} Beals R., Coifmann R., The spectral problem for the Davey-
Stewartson and Ishomori hierarchies, In: "Nonlinear Evolution Equations:
Integrability and Spectral methodes", Proc.Workshop, Como, Italy 1988, Proc.
Nonlinear Sci., 15-23 (1990)
\item{[  BD]} Berenstein C., Dostal M., Some remarks on convolution equations,
 Annales de

\item{      } l'Institut Fourier, {\bf 23}, 55-73 (1973)
\item{[BLMP]} Boiti M., Leon J., Manna M., Pempinelli F., On a spectral
transform of a KDV-like equation related to the Schr\"odinger operator in
the plane, Inverse problems {\bf 3}, 25-36 (1987)
\item{[  Bu]} Bukhgeim A.L., Recovering a potential from the Cauchy data
in the two-dimensional case, J.Inv.Ill-posed Problems, {\bf 16}, (2008)
\item{[   C]} Calderon A.P., On an inverse boundary problem. In: Seminar
on Numerical Analysis and its Applications to Continuum Physics, Soc.
Brasiliera de Matematica, Rio de Janeiro, pp. 61-73 (1980)
\item{[   D]} Druskin V.L., The unique solution of the inverse problem in
electrical surveying and electrical well logging for piecewise-constant
conductivity, Physics of the Solid Earth {\bf 18}(1), 51-53 (1982)
\item{[ DKN]} Dubrovin B.A., Krichever I.M., Novikov S.P.,
The Schr\"odinger equation in a periodic field and Riemann surfaces,
Dokl.Akad.Nauk SSSR {\bf 229}, 15-18 (1976) (in Russian), Sov.Math.Dokl.,
{\bf 17}, 947-951, (1976)
\item{[   E]} Ehrenpreis L., Solutions of some problems of division II,
Amer. J. Math. {\bf 77}, 286-292 (1955)
\item{[  F1]} Faddeev L.D., Increasing solutions of the Schr\"odinger
equation, Dokl.Akad.Nauk SSSR, {\bf 165}, 514-517 (1965)
(in Russian), Sov.Phys.Dokl. {\bf 10}, 1033-1035 (1966)
\item{[  F2]} Faddeev L.D., The inverse problem in the quantum theory of
scattering II, Curr.Probl.
\item{      } Math. {\bf 3}, 93-180 (1974) (in Russian),
J.Sov.Math. {\bf 5}, 334-396 (1976)
\item{[  Fe]} Fedorjuk M.V., Asymptotic: integrals and series, M.Nauka
(1987) (in Russian)
\item{[  Ga]} Garsia A.M., An imbedding of closed Riemann surfaces in
euclidean space, Comm.
\item{      } Math.Helv. {\bf 35}, 93-110 (1961)
\item{[  Ge]} Gelfand I.M., Some problems of functional analysis and algebra,
In: Proc.Int.Congr.
\item{      } Math., Amsterdam, pp.253-276 (1954)
\item{[  GH]} Griffiths Ph., Harris J., Principles of algebraic geometry,
John Wiley, 1978
\item{[  GN]} Grinevich P.G., Novikov S.P., Two-dimensional inverse
scattering problem for negative energies and generalized analytic functions,
Funktsional Anal i Prilozhen. {\bf 22}(1), 23-33 (1988)
\item{[  GT]} Guillarmou C., Tzou L., Calderon inverse problem for the
Schr\"odinger operator on Riemann surfaces, arXiv:0904.3804 (2009) v.1
\item{[  Ha]} Hartshorne R., Algebraic geometry, Springer-Verlag,  (1977)
\item{[  He]} Henkin G.M., Cauchy-Pompeiu type formulas for $\bar\pa$ on
affine algebraic Riemann surfaces and some applications,
arXiv:0804.3761, (2008) v.1 (2010) v.2
\item{[  HP]} Henkin G.M., Polyakov P.L., Homotopy formulas for the $\bar\pa$-
operator on ${\C}P^n$ and the Radon-Penrose transform, Math.USSR
Izvestiya {\bf 28}, 555-587 (1987)
\item{[ HM1]} Henkin G.M., Michel V., On the explicit reconstruction of a
Riemann surface from its Dirichlet-to-Neumann operator, GAFA, Geom.Funct.Anal.
{\bf 17}, 116-155 (2007)
\item{[ HM2]} Henkin G.M., Michel V., Inverse conductivity problem on
Riemann surfaces, J.Geom.
\item{      } Anal. {\bf 18}, 1033-1052 (2008)
\item{[ Ho]} Hodge W., The theory and applications of harmonic integrals,
Cambridge Univ.Press, 1952
\item{[  H\"o]} H\"ormander L., The analysis of linear partial differential
operators I, Springer 1990
\item{[  KV]} Kohn R., Vogelius M., Determining conductivity by boundary
measurements II,
\item{      } Comm.Pure Appl.Math. {\bf 38}, 644-667 (1985)
\item{[  Na]} Nachman A., Global uniqueness for a two-dimensional inverse
boundary problem, Ann. of Math. {\bf 143}, 71-96 (1996)
\item{[  N1]} Novikov R., Reconstruction of a two-dimensional Schr\"odinger
operator from the scattering amplitude at fixed energy,
Funktsional Anal i Prilozhen. {\bf 20}(3), 90-91 (1986) (in Russian)
\item{[  N2]} Novikov R., Multidimensional inverse spectral problem for the
equation $-\Delta\psi+(v(x)-E u(x))\psi=0$,
Funktsional Anal i Prilozhen. {\bf 22}(4), 11-22 (1988) (in Russian)
\item{[  N3]} Novikov R., The inverse scattering problem on a fixed energy
level for the two-dimensional Schr\"odinger operator, J.Funct.Anal.{\bf 103}
(2), 409-463 (1992)
\item{[  NV]} Novikov S.P., Veselov A.P., Two-dimensional Schr\"odinger
operators in periodic fields, Current Problems in Math. {\bf 23}, 3-32 (1983)
(in Russian)
\item{[  Ro]} Rodin Y., Generalized analytic functions on Riemann surfaces,
Lecture Notes Math., {\bf 1288}, Springer (1987)
\item{[  Ru]} R\"uedy R.A., Embeddings of open Riemann surfaces, Comm.Math.
Helv. {\bf 46}, 214-225 (1971)
\item{[  SU]} Sylvester I., Uhlmann G., A uniqueness theorem for an inverse
boundary value problem in electrical prospection, Comm.Pure Appl.Math.
{\bf 39}, 91-112 (1986)
\item{[  Ts]} Tsai T.Y., The Schr\"odinger operator in the plane, Inverse
Problems {\bf 9}, 763-787 (1993)
\item{[   V]} Vekua I.N., Generalized analytic functions, Pergamon, (1962)

\end